\documentclass{article}
\usepackage[latin1]{inputenc}
\title{A compactness result for Landau state in thin-film micromagnetics}
\author{{\Large Radu Ignat}\footnote{Laboratoire de Math\'ematiques, Universit\'e
Paris-Sud 11, b\^at. 425, 91405 Orsay, France \,\,(e-mail:
Radu.Ignat@math.u-psud.fr)} \and {\Large Felix
Otto}\footnote{Max Planck Institute for Mathematics in the Sciences, Inselstr. 22-26,
D-04103 Leipzig,
Germany
 \,\,(e-mail: {otto@mis.mpg.de})}}

\usepackage[mathscr]{euscript} \usepackage{graphicx, pstricks, pst-plot, pst-math} 
\usepackage{amsfonts}
\usepackage{amsmath, latexsym}
\usepackage{amssymb,verbatim}
\usepackage{theorem}
\usepackage{color}

\newtheorem{lem}  {Lemma}
\newtheorem{pro}    {Proposition}
\newtheorem{thm}    {Theorem}

\newtheorem{df}     {Definition}

\theorembodyfont{\rmfamily}
\newtheorem{rem}{Remark}

\newcommand{\bbbb}{B}
\newcommand{\ga}{\gamma}
\newcommand{\dte}{\delta \theta}
\newcommand{\RR}{\mathbb{R}}

\newcommand{\eps}{\varepsilon}

\newcommand{\ds}{\displaystyle}
\newcommand{\h}{{\mathcal{H}}}
\newcommand{\caf}{{\mathcal{F}}}
\newcommand{\proof}[1]{\par\medskip\noindent{\bf Proof#1.}}
\newcommand{\qed}{\hfill$\square$}

\newcommand{\be}{\begin{equation}}
\newcommand{\ee}{\end{equation}}

\newcommand{\ph}{\varphi_\eps}

\newcommand{\nd}{\noindent}
\newcommand{\bu}{B_{1/200}}

\newcommand{\Ix}{I_{x_1}}
\newcommand{\NN}{\mathbb{N}}

\newcommand{\Om}{\Omega}
\newcommand{\dist}{\mathop{\rm dist \,}}
\renewcommand{\t}{\theta_\lambda}
\newcommand{\lam}{\lambda}
\newcommand{\f}{\varphi}

\newcommand{\degr}{\operatorname{deg}}

\newcommand{\vfi}{\varphi}

\newcommand{\ze}{\zeta_\eps}

\newcommand{\tm}{\tilde m'}
\newcommand{\tmp}{\tilde m'(\Phi(x))}
\newcommand{\ta}{{\bf \ddot{w}}(s) }
\newcommand{\na}{{\bf \nu}(s) }
\newcommand{\wa}{{\bf \dot{w}}(s) }
\newcommand{\al}{{\alpha_s(t) }}
\newcommand{\aln}{{\alpha_s(-t) }}
\newcommand{\hd}{\dot{H}^{1/2}(\RR)}
\newcommand{\hu}{\dot{H}^1(\RR^2)}
\begin{document}

\maketitle

\begin{abstract}
 We deal with a nonconvex and nonlocal variational problem coming from thin-film micromagnetics.
 It consists in a free-energy functional depending on two small parameters $\eps$ and $\eta$ and defined
 over vector fields $m : \Omega\subset \RR^2 \to S^2$ that are tangent at the boundary $\partial \Omega$.
 We are interested in the behavior of minimizers as $\eps, \eta \to 0$. They tend to be in-plane away
 from a region of length scale $\eps$ (generically, an interior vortex ball or two boundary vortex balls) and of vanishing divergence, so that $S^1-$transition layers of length scale $\eta$ (N\'eel walls) are
 enforced by the boundary condition. We first prove an upper bound for the minimal energy that corresponds
 to the cost of a vortex and the configuration of N\'eel walls associated to the viscosity solution, so-called Landau state. 
 Our main result concerns the compactness of vector fields $\{m_{\eps, \eta}\}_{\eps, \eta \downarrow 0}$ of energies close to the Landau state
 in the regime where a vortex is energetically more expensive than a N\'eel wall. Our method uses techniques developed for the Ginzburg-Landau type problems for the concentration of energy on vortex balls, together with 
 an approximation argument of $S^2-$vector fields by $S^1-$vector fields away from the vortex balls.
\end{abstract}



\section{Introduction}

In this paper, we investigate a common pattern of the magnetization in thin ferromagnetic films, called Landau state, that corresponds to the global minimizer of the micromagnetic energy
in a certain regime. For that, we focus on a toy problem rather than on the full physical model: 

Let $\Omega\subset \RR^2$ be a bounded simply-connected
domain with a $C^{1,1}$ boundary corresponding to the horizontal section of a ferromagnetic cylinder of small thickness. Due to the thin film geometry, the variations of the magnetization in the thickness direction are strongly penalized. It motivates us to consider magnetizations that are invariant in the out-of-plane variable, i.e.,
$$m=(m_1, m_2, m_3):{\Omega}\to S^2$$ and they are tangent to the
boundary $\partial \Omega$, i.e., \be \label{cond_tan} m'\cdot \nu=0\quad \textrm{
on } \partial \Omega,\ee where $m'=(m_1,m_2)$ is the in-plane
component of the magnetization and $\nu$ is the normal outer unit vector to
$\partial \Omega$. We consider the following micromagnetic energy functional:
$$E_{\eps, \eta}(m)=\int_\Omega |\nabla m|^2\, dx+\frac{1}{\eps^2}\int_\Omega m_3^2\, dx+\frac 1 \eta \int_{\RR^2}
||\nabla|^{-1/2}(\nabla \cdot m')|^2\,
dx,
$$
where $\eps$ and $\eta$ are two small positive parameters (standing for
the size of the vortex core and the N\'eel wall core, respectively). Here,
$x=(x_1,x_2)$ are the in-plane variables with the differential
operator
$$\nabla=(\partial_{x_1}, \partial_{x_2}).$$

The first term of $E_{\eps, \eta}(m)$ stands for the exchange energy. The second term corresponds to the stray-field energy penalizing the top and bottom surface charges $m_3$ of the magnetic cylinder, while the last term counts the stray-field energy penalizing the volume charges $\nabla \cdot m'$ where we will always think of $$m'\equiv m' {\bf 1}_\Omega$$ as
being extended by $0$ outside $\Omega$. For more physical details, we refer to Section \ref{physi}. 

Note that the non-local term in the energy is given by the homogeneous $\dot H^{-1/2}-$seminorm of the
in-plane divergence $\nabla \cdot m'$ that writes in the Fourier
space as: \be \label{defH-12}  \|\nabla \cdot m'
  \|^2_{\dot H^{-1/2}(\RR^2)}=\int_{\RR^2} \bigg|\, |\nabla|^{-1/2}(\nabla
\cdot m')\bigg|^2\, dx:=\int_{\RR^2} \frac{1}{|\xi|} |{\cal F}(\nabla
\cdot m')|^2\, d\xi. \ee Also observe that the boundary condition
\eqref{cond_tan} is necessary so that \eqref{defH-12} is finite
since
$$\nabla \cdot m'=(\nabla \cdot m'){\bf 1}_\Omega+ (m'\cdot
\nu) {\bf 1}_{\partial \Omega}\quad \textrm{ in } \RR^2$$(see Proposition \ref{prop_app} in Appendix). 

We are interested in the asymptotic behavior of minimizers of the
energy $E_{\eps, \eta}$ in the regime $$\eps\ll 1\quad \textrm{ and }\quad \eta\ll 1.$$ 
The main features of this variational model resides in the nonconvex constraint on
the magnetization $|m|=1$ and the nonlocality of the stray-field interaction. The
competition of these effects with the quantum mechanical exchange effect leads to a rich pattern formation for the stable
states of the magnetization. Generically, a pattern of a stable state consists in
large uniformly magnetized regions ({\it magnetic domains}) separated by narrow
smooth transition layers ({\it wall domains}) where the magnetization varies rapidly.
The characteristic wall domains observed in thin ferromagnetic films are the N\'eel walls (corresponding to a one-dimensional in-plane
rotation connecting two directions of the magnetization) together with topological defects standing for interior vortices (called Bloch lines)
and micromagnetic boundary vortices. 

The existence of line singularities at the mesoscopic level of the magnetization in thin films can be explained by the principle of pole avoidance.
For this
discussion, we first neglect the exchange term in $E_{\eps, \eta}$. The stray-field energy
will try to enforce in-plane configurations, i.e., $m_3=0$ in $\Omega$, together with the divergence-free condition for $m'$, i.e.,
$\nabla\cdot m'=0$ in $\Omega$. Together with \eqref{cond_tan}, we arrive at \be
\label{eikonal1} |m'|=1, \, \nabla \cdot m'=0 \,\textrm{
in }\, \Omega \, \textrm{ and } \, m'\cdot \nu=0\quad \textrm{
on } \partial \Omega.\ee We notice that the conditions in
\eqref{eikonal1} are too rigid for smooth magnetization $m'$.
This can be seen by writing $m'=\nabla^\perp\psi$ with the help
of a ``stream function" $\psi$. Then up to an additive constant, \eqref{eikonal1} implies that
$\psi$ is a solution of the Dirichlet problem for the eikonal
equation: \be \label{eikona} |\nabla\psi|=1 \,\textrm{ in }\,
\Omega \quad  \textrm{
and }  \psi=0 \textrm{
on } \partial \Omega.\ee The method of characteristics yields the nonexistence
of smooth solutions of \eqref{eikona}. But there are many continuous solutions that satisfy
\eqref{eikona} away from a set of vanishing Lebesgue measure. One
of them is the ``viscosity solution" given by the distance
function
$$\psi(x)=\dist(x, \partial \Omega')$$ that corresponds to the so-called
Landau state for the magnetization $m'$. Hence, the boundary conditions \eqref{cond_tan} are expected to induce
line-singularities for solutions $m'$ that are an idealization of wall domains at the mesoscopic level. At the microscopic
level, they are replaced by smooth transition layers, the N\'eel walls, where the
magnetization varies very quickly on a small length scale $\eta$. Note that the normal component
of $m'$ does not jump across these discontinuity lines (because of
\eqref{eikonal1}); therefore, the normal vector of the mesoscopic wall
is determined by the angle between the mesoscopic levels of the magnetization in the
adjacent domains (called angle wall). Now, taking into account the contribution of the exchange effect, 
the energy scaling
per unit length of a N\'eel wall of angle $2\theta$ (with $\theta\in(0, \frac \pi 2])$) is given in DeSimone, Kohn, M{\"u}ller\& Otto \cite{DKMObook}, Ignat \&Otto \cite{IO} (see also Ignat \cite{IGNAT_gamma}):
\be
\label{ene_nee}
\frac{\pi (1-\cos \theta)^2+o(1)}{\eta |\log \eta|}\, \quad \textrm { as } \eta \to 0. \ee

The formation of interior or boundary vortices is explained by the competition between the exchange energy and the penalization of the 
$m_3-$component for configurations tangent at the boundary. Indeed, there is no $S^1-$configuration that is of finite exchange energy and satisfies \eqref{cond_tan}. There are only two possible 
situations: 
If $m'$ does not vanish on $\partial \Omega$, than \eqref{cond_tan} implies that $m'$ carries a nonzero topological degree, $\degr(m', \partial \Omega)=\pm 1$. In this case,
we expect the nucleation of an interior vortex of core-scale $\eps$. The scaling of the vortex energy is related to the
minimal Ginzburg-Landau (GL) energy (see Bethuel, Brezis \& Helein \cite{BBH2}):
\be
\label{en_GL_in}
\mathop{\min_{m'\in H^1(\Omega, \RR^2)}}_{m'=\nu^\perp \textrm{ on } \partial
\Omega} \int_\Omega g_\eps(m')\,
dx=(2\pi+o(1)) |\log \eps| \quad \textrm { as } \eps \to 0,\ee where the GL density energy is
given in the following: \be \label{GLdef} g_\eps(m')=|\nabla
m'|^2+\frac{1}{\eps^2}\left(1-|m'|^2\right)^2.\ee
(Here, we denote $\nu^\perp=(-\nu_2, \nu_1)$.)
The second situation consists in having zeros of $m'$ on the boundary. Therefore, we expect that boundary vortices do appear. Roughly speaking, they correspond to "half" of an interior vortex 
where the vector field $m'$ is tangent at the boundary; therefore they are different from the micromagnetic
boundary vortices analyzed by
Kurzke\cite{Kurzke-2006} and Moser \cite{Moser} (see details in Section \ref{physi}).
Remark the importance of the regularity of $\partial \Omega$ in estimate \eqref{en_GL_in}. In fact, if $\partial \Omega$ has a corner and the boundary condition
$m'=\nu^\perp$ on $\partial
\Omega$ in \eqref{en_GL_in} is relaxed to \eqref{cond_tan}, then estimate \eqref{en_GL_in} does not hold anymore, it depends on the angle of the
corner (see Proposition \ref{exec1} and Remark \ref{remu}).
Therefore, at the microscopic level, topological point defects do appear in the Landau state pattern and are induced by \eqref{cond_tan}. 

The aim of the paper is to show compactness of magnetizations of energy $E_{\eps, \eta}$ close to the Landau state in order to rigorously justify the
limit behavior \eqref{eikonal1}: the delicate issue consists in having the constraint $|m|=1$ conserved in the limit.  For that, we have to evaluate the energetic cost of the Landau state. We expect that the leading order energy of a Landau state 
is given by
the topological point defects and N\'eel walls. 
The Landau state configuration consists in several N\'eel walls and either one interior Bloch line
or two "half" Bloch lines placed at the boundary of the sample $\Omega$. Therefore, by \eqref{ene_nee} and \eqref{en_GL_in}, we expect that
the energy of the Landau state has the following order:
\be
\label{upp_conj}
2\pi |\log
\eps|+\frac{A}{\eta |\log \eta|},\ee for some positive $A>0$ depending on the length and angle of N\'eel walls.

\section{Main results}

\begin{figure}[htbp]
\center
\includegraphics[scale=0.3,
width=0.3\textwidth]{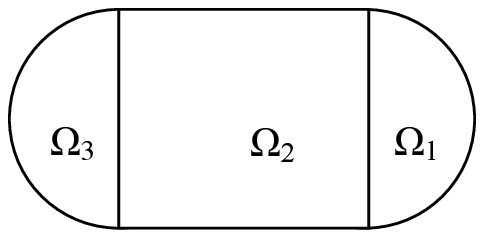} \caption{Stadium} \label{stadu}
\end{figure}

First of all, we want to rigorously prove the upper bound \eqref{upp_conj} for the Landau state. Our result gives the exact leading order energy of the Landau state in
the case of a domain $\Omega$ of a "stadium" shape (see Figure \ref{stadu}). 
Note that the Landau state of a stadium consists in a single N\'eel wall of
$180^\circ$ (in our example, the length of the wall is equal to $2$, so that $A=2\pi$ in \eqref{upp_conj}). 
\begin{thm}
\label{thm_upper} Let $\Omega=\Omega_1\cup \Omega_2\cup \Omega_3$
be the following "stadium" shape domain:
\begin{align*}
\Omega_1&=\{x=(x_1, x_2)\in \RR^2\, :\, |x-(1,0)|<1, x_1\geq
1\},\\
\Omega_2&=(-1,1)\times (-1, 1),\\
\Omega_3&=\{x=(x_1, x_2)\in \RR^2\, :\, |x-(-1,0)|<1, x_1\leq
-1\}.
\end{align*} In the regime $\eps \ll  \eta \ll 1$, there exists a $C^1$ vector field
$m_{\eps,\eta}:{\Omega}\to {S^2}$ that satisfies \eqref{cond_tan}
and
\be
\label{estime1}
E_{\eps, \eta}(m_{\eps,\eta})\leq 2\pi |\log
\eps|+\frac{2\pi+o(1)}{\eta |\log \eta|}
\quad \textrm{ as } \quad \eta\downarrow
0.\ee
\end{thm}

\bigskip

Observe that the vortex energy in the above estimate is relevant only if a vortex costs at least as much as a N\'eel wall, i.e., $\frac{1}{\eta |\log \eta|} \lesssim |\log
\eps|$ (otherwise, the vortex energy would be absorbed by the term $o\big(\frac{1}{\eta |\log \eta|}\big)$\,). This regimes leads to a size $\eps$ of the vortex core exponentially smaller than the size of the N\'eel wall core $\eta$ (see Remark \ref{rmk1}).

\bigskip

{\noindent \bf Notation:} We always denote $a\ll b$ if $\frac a b \to 0$; also, $a\lesssim b$ if $a\leq Cb$ for some universal constant $C>0$.

\bigskip

Now we state our main result on the compactness of the $S^2-$valued magnetizations that
have energies near the Landau state. The issue consists in rigorously justifying that the constraint $|m|=1$ is conserved by the limit configurations as $\eps, \eta \to 0$. 
The regime where we prove our result corresponds to the case where a topological defect is energetically more expensive than the N\'eel wall, that is coherent with the regime where \eqref{estime1} holds.

\begin{thm}
\label{thm_main} Let $\alpha\in (0, \frac 1 2)$ be an
arbitrary constant. We consider the following regime between the small parameters $\eps, \eta \ll 1$: 
\be \label{regim_nou}
\eps^{1/2}\lesssim \eta,\ee
\be \label{regim_cond11} \log |\log \eps|\lesssim \frac{1}{\eta |\log
\eta |}.\ee
For each $\eps$ and $\eta$, we consider $C^1$ vector fields
$m_{\eps, \eta}:{\Omega}\to S^2$ that satisfy \eqref{cond_tan}
and 

\begin{minipage}{0.55\linewidth}
\hspace{1cm}
\quad $E_{\eps, \eta}(m_{\eps, \eta})-2\pi |\log \eps| \quad \quad$ {\Huge\{} 
\end{minipage}
\hfill
\begin{minipage}{0.35\linewidth} 
 \hspace{1cm} \begin{align} \label{en_level}  &\leq  2\pi \alpha |\log
\eps|\\ \label{en_level2} &\lesssim \frac{1}{\eta |\log \eta|}
 \end{align}
\end{minipage}
 
\noindent Then the family
$\{m_{\eps, \eta}\}_{\eps, \eta\downarrow 0}$ is relatively compact in
$L^1(\Omega, S^2)$ and any accumulation point $m:\Omega\to S^2$
satisfies \be \label{limit_prop} m_3=0, \,\,|m'|=1 \textrm{ a.e.
in $\Omega \quad$ and$ \quad \nabla \cdot m'=0$ distributionally
in $\RR^2$}.\ee
\end{thm}

\bigskip

The proof of compactness is based on an argument of approximating
$S^2-$valued vector fields by $S^1$-valued vector fields away from a
small defect region. This small region consists in either one interior vortex or two boundary vortices. The detection of this region is done in Theorem
\ref{detect_vort} and uses some topological methods due to Jerrard
\cite{Jer99} and Sandier \cite{Sa} for the concentration of the
Ginzburg-Landau energy around vortices (see also Lin \cite{Lin}, Sandier \& Serfaty \cite{SSbook} ). Away from this small region,
the energy level {\it only} allows for line singularities.
Therefore, the compactness result for $S^1-$valued vector fields in
\cite{IO} applies.

Let us discuss the assumptions \eqref{regim_nou}, \eqref{regim_cond11},  \eqref{en_level} \&  \eqref{en_level2}. Inequality \eqref{en_level2} assures that cutting out the topological defect (one vortex or two boundary vortices), the remaining energy rescaled at the energetic level of N\'eel walls is uniformly bounded. Inequality \eqref{en_level} together with the choice of $\alpha<\frac 1 2$ mean that the energy cannot support three "half" interior vortices and is precisely explained in Theorem \ref{detect_vort} below.
Inequality \eqref{regim_cond11} is imposed due to our method to detect a boundary vortex: it leads to a loss of energy of order $O(\log |\log \eps|)$ with respect to the expected half energy of a interior vortex, i.e., $\pi |\log \eps|$ (see Theorem \ref{detect_vort} and Proposition \ref{exec1}). This amount of energy could leave room for configurations of N\'eel walls that may distroy the compactness of $|m'|=1$. Therefore, to avoid this scenario, \eqref{regim_cond11} is imposed. The regime \eqref{regim_nou} is rather technical: it is needed in the approximation argument of $S^2-$valued vector fields by $S^1-$valued vector fields away from the vortex balls. In fact, starting from the values of $m'$ on a square grid of size $\eps^\beta$, the approximation argument requires zero degree of $m'$ on each cell, leading to the condition $\beta<1-\alpha$ (see Lemma \ref{lemcell}); furthermore, the condition $\eps^\beta\lesssim \eta$ is needed in order that the approximating $S^1-$valued vector fields induce a stray field energy of the same order of $m'$ (see \eqref{star2}). Therefore, \eqref{regim_nou}
can be improved to a larger regime $$\eps^\beta\lesssim \eta \quad \textrm{ for any } \quad \beta<1-\alpha$$ as presented in the proof (Theorem \ref{thm_main} is stated for the value $\beta=1/2$ which is the universal choice for every $\alpha<1/2$). However, this slightly improved condition is weaker than the complete regime implied by 
 \eqref{en_level} as explained in the following remark.
 
 \begin{rem}
 \label{rmk1}
 Any limit configuration $m'$ satisfies \eqref{limit_prop}. If $\Om$ is a bounded simply-connected domain different than discs, $m'$ has at least one ridge (line-singularity) that corresponds
 to a N\'eel wall. Therefore, the minimal energy verifies $\mathop{\min}_{\eqref{cond_tan}} E_{\eps, \eta}-2\pi |\log \eps|\gtrsim \frac{1}{\eta |\log \eta|}$. Combining with \eqref{en_level}, it follows that
 $$\frac{1}{\eta |\log \eta|} \lesssim |\log
\eps|;$$
in particular, $\eps\lesssim e^{-\frac{1}{\eta |\log \eta|}}$, i.e., the core of the vortex is exponentially smaller than the core of the N\'eel wall. However, in the proof of Theorem \ref{thm_main}, this much stronger constraint with respect to \eqref{regim_nou} is not needed. 
 \end{rem}

We prove the following result of the concentration of Ginzburg-Landau energy around one interior vortex or two boundary vortices for vector fields tangent at the boundary:

\begin{thm}
\label{detect_vort} Let $\alpha\in (0, \frac 1 2)$ and $\Omega\subset \RR^2$ be a bounded simply-connected
domain with a $C^{1,1}$ boundary. There exists $\eps_0=\eps_0(\alpha, \partial \Omega)>0$ such that
for every $0<\eps<\eps_0$, if $m':{\Omega}\to \overline{B^2}$ is a $C^1$ vector field that satisfies \eqref{cond_tan}
and \be \label{en_vort} \int_\Omega g_\eps(m')\, dx\leq 2\pi
(1+\alpha)|\log \eps|,\ee then there exists either a ball $B(x^*_1,
r^*)\subset \Omega$ (called vortex ball) with $r^*=\frac{1}{|\log \eps|^3}$ and
\be
\label{estvoin}
\noindent \int_{B(x^*_1, r^*)}g_\eps(m')\, dx\geq 2\pi |\log \frac{r^*}{\eps}|-C,\ee or two balls $B(x^*_2, r^*)$ and  $B(x^*_3, r^*)$ (called boundary vortex balls) with $x^*_2, x^*_3\in \partial \Omega$ and
\be
\label{estvobd}
\int_{(B(x^*_2, r^*)\cup B(x^*_3, r^*))\cap \Omega}g_\eps(m')\, dx\geq 2\pi |\log \frac{r^*}{\eps}|-C,\ee
where $C=C(\alpha, \partial \Omega)>0$ is a constant depending only on $\alpha$ and on the geometry of $\partial \Omega$.
\end{thm}

The condition $\alpha<1/2$ is needed in our proof. In fact, if no topological defect exists in the interior (in which case, condition \eqref{cond_tan} induces boundary vortices), we perform a
mirror-reflection extension of $m'$ outside the domain. Roughly speaking, the GL energy in the extended domain doubles, i.e.,
it is of order $2\pi (2+2\alpha)|\log \eps|$ and the degree at the new boundary is equal to two; in order to avoid the formation of three interior vortices in the extended region, we should impose $2+2\alpha<3$,  i.e.,
$\alpha<1/2$.

Notice that the Ginzburg-Landau energy concentration for a boundary vortex in \eqref{estvobd} has a cost of order $\pi |\log \eps|-C\log |\log \eps|$ provided that the boundary has regularity
$C^{1,1}$. We conjecture that the same energetic cost for a boundary vortex holds true if the boundary has regularity $C^{1, \beta}$, $\beta\in (0,1)$. However, if the boundary regularity is only $C^1$, then the energetic cost of a
boundary vortex may decrease to $(\pi-\frac{C}{\log |\log \eps|})|\log \eps|$ where $C>0$ is a universal constant.
This indicates that the loss of energy of order $\log |\log \eps|$ in \eqref{estvobd} could occur for
boundary vortices for $C^{1,\beta}$ boundary regularity and the order of this loss increases to $\frac{|\log \eps|}{\log |\log \eps|}$ for $C^1$ boundaries as $\beta\to 0$. This claim is 
supported by the following example for a $C^1$ boundary domain:
\begin{pro}
\label{exec1}
We consider in polar coordinates the following $C^1$ domain $\Omega=\{(r, \theta)\,:\, r\in (0, \frac{1}{20}), \, \,  |\theta|<\ga(r)=\frac \pi 2-\frac{1}{\log \log \frac{1}{r}}\}$.
For every $0<\eps<1$, there exists a $C^1-$function $m'_\eps:\Omega\cap \bu\to \RR^2$ that satisfies \eqref{cond_tan} on $\partial \Omega\cap \bu$ and
$$\int_{\Omega\cap \bu} g_\eps(m'_\eps)\, dx\leq (\pi-\frac{C}{\log |\log \eps|}) |\log \eps|,$$
where $C>0$ is some universal positive constant (independent of $\eps$).
\end{pro}

The outline of the paper is as follows. In Section \ref{physi}, we present the physical context of our toy problem. In the next section, we recall two results that we need for the proof
of our results:
a compactness result for $S^1-$valued magnetizations and the concentration of the Ginzburg-Landau energy on vortex balls. In Section \ref{section5}, we prove Theorem \ref{detect_vort} and Proposition \ref{exec1}. In Section \ref{section6}, we give the proof of our main result in Theorem \ref{thm_main}. In Section \ref{section_up}, we show the upper bound for the stadium domain stated in Theorem
\ref{thm_upper}.
In Appendix, we prove that \eqref{cond_tan} is a necessary 
condition for our configurations to have a finite stray field energy.

\section{Physical context}
\label{physi}

In this section we explain the physical context of our model in thin-film micromagnetics. We consider 
a ferromagnetic sample of cylinder shape,
  i.e.  $$\omega = \omega' \times (0,t)$$ where $\omega' \subset \RR^2$ is the section of the magnetic sample of length $\ell$ and $t$ is the thickness of the cylinder. 
The microscopic behavior of the magnetic body is described by a three-dimensional unit-length vector field $m=(m', m_3):\omega \to S^2$, called magnetization.
The observed ground state of the magnetization is a minimizer of 
the micromagnetic energy that we write here in the absence of anisotropy and external magnetic field:
\begin{align} \label{e-3d} %
  E^{3d}(m) = d^2 \ \int_{\omega} |(\nabla, \frac{\partial}{\partial z}) m|^2 \ dx dz
  + \int_{\RR^3} | (\nabla, \frac{\partial}{\partial z})U(m)|^2 \ dx dz.
 \end{align}
 The parameter $d$ of the material is called
 exchange length and is of order of nanometers. The stray-field potential
 $U(m):\RR^3\to \RR$ is defined by static Maxwell's equation in the weak sense:
  \begin{align} \label{strayfield-3d} %
    \int_{\RR^3} (\nabla, \frac{\partial}{\partial z}) U(m) \cdot (\nabla,
    \frac{\partial}{\partial z}) \zeta \ dx dz = \int_{\RR^3} (\nabla,
    \frac{\partial}{\partial z}) \cdot (m {\bf 1}_\omega) \ \zeta \ dx dz, \quad \textrm{ for every }
    \zeta\in C^\infty_c(\RR^3).
  \end{align}
  Instead of the three length scales  $\ell$, $t$ and $d$ of the physical model, we introduce two dimensionless parameters:
  \begin{align*}
   \eps := \frac d \ell && \text{and} && \eta:= \frac{d^2}{\ell t} 
  \end{align*}
 (standing for the size of the core of the Bloch line and the N\'eel wall, respectively).  
 
 \medskip
 
 {\nd \bf Thin-film reduction.}  We consider the thin-film approximation of the full energy \eqref{e-3d} in the following regime:
  \begin{align} \label{regime_thin} %
    \eps \ll \eta \ll 1
  \end{align}
(equivalently,  $t \ll d \ll \ell$). The assumption $t\ll d$ implies that in-plane transitions (N\'eel walls) are preferred to out-of-plane transitions (asymmetric Bloch walls)  between two mesoscopic directions of the magnetization (see Otto \cite{Otto_cross_over}). The hypothesis $d \ll \ell$ assures that constant configurations in general are not global minimizers (see DeSimone \cite{DeSimone_mec}). 

The
 main issue is the asymptotic behavior of the energy in the regime of thin
 films.  We first nondimensionalize in length with respect to $\ell$, i.e. $(\bar x,
 \bar z)=(\frac{x}{\ell}, \frac{z}{\ell})$, $\Omega=\frac{\omega'}{\ell}$, $\bar m(\bar x, \bar z)=m(x,z)$, $\bar
 U(\bar m)(\bar x, \bar z)= \frac 1 \ell U(m)(x,z)$ and then we renormalize the energy $\bar E^{3d}(\bar
 m)=\frac{1}{d^2 t} E^{3d}(m)$.  Omitting the $\ \bar{} \ $, we get
 \begin{align} \label{eee}
   E^{3d}(m) = \frac{\eta}{\eps^2} \ \int_{\Omega \times (0, \frac{\eps^2}{\eta})} |(\nabla,
   \frac{\partial}{\partial z}) m|^2 \ dx dz + \frac{\eta}{\eps^4}\int_{\RR^3} | (\nabla,
   \frac{\partial}{\partial z})U(m)|^2 \ dx dz.
 \end{align}
In the regime \eqref{regime_thin}, the penalization of exchange energy 
 enforces the following constraints for the minimizers:  
 \begin{enumerate}
 \item[(a)] $m$ varies on length scales $\gg \frac{\eps^2}{\eta}$.
 \item[(b)] $m = m(x)$, i.e. $m$ is $z-$invariant.
 \end{enumerate}
 With these assumptions, \eqref{eee} can be approximated by the following reduced
 energy $E^{red}$ (see DeSimone, Kohn, M\"uller \& Otto \cite{DKMO1}, Kohn \& Slastikov \cite{KohnSlastikov-2005}):
 \begin{align}  %
     \nonumber
   E^{red}(m) &= \int_{\Omega} |\nabla m|^2 \ dx\\
  \label{enmea}
 & +\frac{1}{\eps^2} \int_{\Omega} m_3^2\, dx
   + \frac{|\log \frac{\eps^2}{\eta}|}{2\pi \eta} \ \int_{\partial \Omega} (m' \cdot \nu)^2 \ d\h^1
  +\ \frac{1}{2\eta}
   \|(\nabla \cdot m')_{ac} \|^2_{\dot H^{-1/2}(\RR^2)}.     
   \end{align}
 The above formula follows by solving the stray field equation \eqref{strayfield-3d} in the regime \eqref{regime_thin}: indeed, for $z-$invariant configurations $m$, the Fourier transform
in the in-plane variables $x=(x_1, x_2)$ turns \eqref{strayfield-3d} into a second order ODE in the $z-$variable that can be solved explicitly (see
 \cite{KohnSlastikov-2005}, \cite{IGNAT-XEDP}). Then, due to the above
 assumption a) and to the regime \eqref{regime_thin}, the stray-field energy asymptotically
 decomposes into three terms as written in \eqref{enmea}: the first term in \eqref{enmea} is penalizing the surface charges $m_3$ on the top and bottom of the
 cylinder, a second
 term counts the lateral charges $m' \cdot \nu$ in the $L^2-$norm, as well as the
 third term that penalizes the volume
 charges $(\nabla \cdot m')_{ac}:=(\nabla \cdot m') {\bf 1_\Omega}$ as a homogeneous $\dot H^{-1/2}-$seminorm. In fact, the last term corresponds to
 the stray-field energy created by a three-dimensional vector field $h_{ac}(m)$ defined as
 \begin{equation*}
  h_{ac}(m)=(\nabla, \frac{\partial}{\partial z}) U_{ac}(m) :\RR^3\to \RR^3, 
\end{equation*}
that satisfies:
\begin{align*}
  \int_{\RR^3} (\nabla, \frac{\partial}{\partial z}) U_{ac}(m) \cdot (\nabla,
  \frac{\partial}{\partial z}) \zeta \ dx dz = \int_{\RR^2} (\nabla \cdot m')_{ac} \ \zeta \
  dx, \textrm{ for all } \zeta \in
C_c^\infty(\RR^3). 
\end{align*}
Then one has
\begin{align}
\label{eqh2}
  \int_{\RR^3} |h_{ac}(m)|^2 \, dx dz
  =\frac 1 2 \|(\nabla \cdot m)_{ac}
  \|^2_{\dot H^{-1/2}(\RR^2)}.
\end{align}
Note that if \eqref{cond_tan} holds (i.e., no lateral surface charges), then $(\nabla \cdot m)_{ac}=\nabla \cdot (m {\bf 1_\Omega})$ and therefore, $h_{ac}(m)$ induces the stray field energy 
given by \eqref{defH-12}. In fact, \eqref{defH-12} corresponds to the minimal stray field energy in thin films.
More precisely, a stray field
$h=(h',h_3)=(h_1,h_2,h_3):\RR^3\to \RR^3$ is related to the
magnetization $m:\Omega\to S^2$ via the following variational formulation: \be
\label{1c} \int_{\RR^2\times\RR} \left(h'\cdot \nabla
\zeta+h_3\frac{\partial \zeta}{\partial z}\right)\,
dxdz=\int_{\RR^2} \zeta \, \nabla \cdot m' \, dx,\quad \forall
\zeta\in C^\infty_c(\RR^3),\ee where $z$ denotes the out-of-plane
variable in the space $\RR^3$. (As before, $m'\equiv m' {\bf 1}_\Omega$ and $m$ satisfies \eqref{cond_tan}.) Classically, this is,
$$
\begin{cases}
\nabla\cdot h'+\frac{\partial h_3}{\partial z}=0 &\quad \textrm{in }\,\,\RR^3\setminus (\RR^2\times \{0\}),\\
[h_3]=-\nabla\cdot m'&\quad \textrm{on }\,\,\RR^2\times \{0\},
\end{cases}$$
where $[h_3]$ denotes the jump of the out-of-plane component of
$h$ across the horizontal plane $\RR^2\times \{0\}$. Then
\eqref{defH-12} can be expressed as:
$$ \int_{\RR^2}
\left|\,|\nabla|^{-1/2}(\nabla\cdot m')\right|^2\, dx=2\min_{h \, {\rm with }\, \eqref{1c} } \int_{\RR^2\times \RR}|h|^2 \,dxdz.$$
Therefore, $h_{ac}(m)$ is a minimizing stray-field (of vanishing curl) associated with the stray field potential $U_{ac}(m)$.

In our
 regime \eqref{regime_thin}, there are three different structures that typically appear: N\'eel walls, Bloch lines and micromagnetic boundary vortices. We explain
 these structures in the following and compare their respective energies. As we already mentioned,
 a fourth structure, the asymmetric Bloch wall, can appear in thicker films but we do not discuss it here since
 the asymmetric Bloch wall is more expensive than a N\'eel wall 
 if $t \ll d$.

 \medskip

\begin{figure}[htbp]
\center
\includegraphics[scale=0.5,
width=0.5\textwidth]{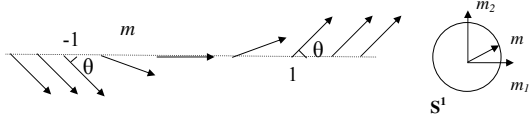} \caption{N\'eel wall
of angle $2\theta$ confined in $[-1,1]$.} \label{Neelwallmod2}
\end{figure}

 {\nd \bf N\'eel walls.} 
The N\'eel wall is a dominant transition layer in thin
ferromagnetic films. It is characterized by a one-dimensional
in-plane rotation connecting two (opposite) directions of the
magnetization. It has two length scales: a small core with fast
varying rotation and two logarithmically decaying tails. In order
for the N\'eel wall to exist, the tails are to be contained and we consider here the confining mechanism of the steric interaction with the sample
edges. Typically, one may consider wall transitions of the form:
$$m=(m_1, m_2):\RR \to S^1 \textrm{ and } m(\pm
t)=\left(
\begin{array}{c}
\cos \theta \vspace{0.2cm}\\ \pm \sin \theta
\end{array}
\right) \textrm{ for } \pm t\geq 1,$$ with $\theta \in [0,\frac \pi 2)$
(see Figure \ref{Neelwallmod2}), whereas the reduced energy functional is:
$$E^{red}(m)=\int_{\RR} |\frac{dm}{dx_1} |^2\, dx_1+
\frac{1}{2\eta} \int_{\RR} \bigg|\, \bigg|\frac{d}{dx_1}\bigg|^{1/2}m_1 \bigg|^2 \, dx_1.$$ 
As $\eta\to 0$, the scale of the N\'eel core is given by $|x_1|\lesssim w_{core}=O(\eta)$ while the two logarithmic
 decaying tails scale as $w_{core}\lesssim|x_1|\lesssim w_{tail}=O(1)$. The energetic cost (by unit length)
 of a N\'eel wall is given by
 \begin{align*}
   E^{red}(\text{N\'eel wall})=O(\frac{1}{\eta |\log \eta|})
 \end{align*}
 with the exact prefactor $\pi(1-\cos \theta)^2/2$ where $2\theta$ is the wall angle (see e.g. \cite{IGNAT_gamma}).

 \medskip

 {\nd \bf Bloch line.} 
 A Bloch line is a regularization of a vortex on the microscopic level of the
 magnetization that becomes out-of-plane at the center. The prototype of a Bloch
 line is given by a vector field
 $$m: B^2\to S^2$$ defined in a circular cross-section $\Omega=B^2$ of
 a thin film and satisfying:
 \begin{align} \label{bloc} %
   \nabla\cdot m'=0 \textrm{ in } B^2 && \text{ and } && m'(x)={x^\perp} \textrm{ on
   }\partial B^2.
 \end{align}
  \begin{figure}[htbp]
  \center  
  \includegraphics[scale=0.5,width=0.5\textwidth]{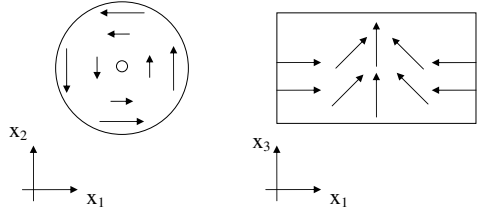}\\
  \caption{Bloch line.} \label{Blochli}
\end{figure}
(For the Bloch line in a thin cylinder, the magnetization is assumed to be invariant in the thickness direction of the film and the
 word ``line'' refers to the vertical direction.) Since the magnetization turns in-plane at
 the boundary of the disk $B^2$ (so, $\degr(m', \partial \Omega)=1$), a localized region is created, that is the core of the Bloch
 line of size $\eps$, where the magnetization becomes perpendicular to the horizontal plane (see Figure
 \ref{Blochli}).
The reduced energy \eqref{enmea} for a configuration \eqref{bloc} writes as: $$
E^{red}(m)=\int_{B^2} |\nabla m|^2 \, dx +\frac{1}{\eps^2} \int_{B^2} m^2_3\, dx.$$
The Bloch line corresponds to the minimizer of this energy under the constraint
\eqref{bloc}. Remark that the reduced energy $E^{red}$ controls the Ginzburg-Landau energy, i.e., 
$$\int_{B^2} g_\eps(m')\, dx\leq  E^{red}(m)$$ since
$|\nabla m'|^2\leq |\nabla (m', m_3)|^2$ and $(1-|m'|^2)^2=m_3^4\leq m_3^2$.
Due to the similarity with
the Ginzburg-Landau type functional, the Bloch line corresponds to the Ginzburg-Landau vortex and the energetic cost of a Bloch line (per unit-length) is given by \eqref{en_GL_in}:
\begin{align*}
  E^{red}(\text{Bloch line})=O({|\log \eps |})
\end{align*}
with the exact prefactor $2\pi$ (see e.g. \cite{IGNAT-XEDP}).

\medskip

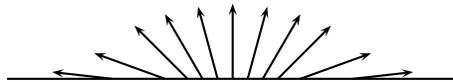
\begin{figure}[ht] %
  \centering
  \begin{pspicture}(-3,0)(3,1) %
    \psset{unit=1cm}
    \psline(-3,0)(3,0)
    \psline{->}(0,0)(0, 1)
    \psline{->}(0.2,0)(0.46, 0.96) 
    \psline{->}(0.4,0)(0.9, 0.86)     
    \psline{->}(0.6,0)(1.307, 0.707) 
    \psline{->}(0.9,0)(1.84, 0.34) 
    \psline{->}(1.6,0)(2.4, 0.09) 
    \psline{->}(-0.2,0)(-0.46, 0.96) 
    \psline{->}(-0.4,0)(-0.9, 0.86)     
    \psline{->}(-0.6,0)(-1.307, 0.707) 
    \psline{->}(-0.9,0)(-1.84, 0.34) 
    \psline{->}(-1.6,0)(-2.4, 0.09) 
    \end{pspicture}
  \caption{A micromagnetic boundary vortex}
  \label{fig-bdrvortex}
\end{figure}

 {\nd \bf Micromagnetic boundary vortex.} 
Next we address micromagnetic boundary vortices. A micromagnetic boundary vortex corresponds to an in-plane
transition of the magnetization along the boundary from $\nu^\perp$ to $-\nu^\perp$,
see Figure \ref{fig-bdrvortex}. The corresponding minimization problem is given by
\begin{align*}
  E^{red}(m) =   \int_{\Omega} |\nabla m|^2 \ dx + \frac{|\log \frac{\eps^2}{\eta}|}{2\pi \eta} \ \int_{\partial \Omega} (m' \cdot \nu)^2 \ d\h^1
\end{align*}
within the set of in--plane magnetizations $m:\Omega\to S^1$. The minimizer of this energy is an harmonic vector field with values in $S^1$
driven by a pair of boundary vortices. These have been analyzed in
\cite{Kurzke-2006, Moser}.  The transition is regularized on the length scale of the exchange part of
the energy, i.e. the core of the boundary vortex has length of size $\frac{\eta}{|\log \frac{\eps^2}{\eta}|}$. The cost of such a transition is given by $$E^{red}(\text{Micromagnetic boundary
  vortex})=O(\bigg|\log \frac{\eta}{|\log \frac{\eps^2}{\eta}|}\bigg|)$$ with exact prefactor $\pi$. (Note that the boundary vortices in Theorem \ref{detect_vort} correspond in fact to "half" Bloch lines where the vector field is tangent at the boundary, i.e., $m' \cdot \nu=0$ on $\partial \Omega$; therefore, their structure is different from the one of micromagnetic boundary vortices, but with the same energetic cost.)  

\medskip

{\nd \bf Claim}: In the regime \eqref{regime_thin}, then
\begin{align*}
  \text{either $\, \quad E^{red}$(Micromagnetic boundary vortex)} &\lesssim \text{$E^{red}$(N\'eel wall)$\,$ }\\ \text{or $\,\quad  E^{red}$(Micromagnetic boundary vortex)}&\lesssim
  \text{$E^{red}$(Bloch line)}.
\end{align*}
Indeed, assume by contradiction that the above statement fails. Then one has
\be
\label{comp11}
\frac{1}{\eta|\log \eta|}\lesssim \bigg|\log \frac{\eta}{|\log \frac{\eps^2}{\eta}|}\bigg|\ee and  
\be
\label{comp22}
\log \frac 1 \eps \lesssim \bigg|\log \frac{\eta}{|\log \frac{\eps^2}{\eta}|}\bigg|.\ee
In the regime \eqref{regime_thin}, one has $\eps^2\ll \eps \ll \eta$, therefore \eqref{comp11} turns into
$$\frac{1}{\eta|\log \eta|} \lesssim \log \log \frac 1 \eps,$$ 
while \eqref{comp22} implies that
$$\log \frac 1 \eps \lesssim \log \frac 1 \eta.$$ 
Now it is easy to see the incompatibility between the last two inequalities as $\eps, \eta \to 0$.

\medskip

{\nd \bf Our toy problem}: The model we presented in the introduction consists in considering configurations without lateral surface charges, i.e., \eqref{cond_tan} holds true.
In this case, our energy functional $E_{\eps, 2\eta}(m)$ coincides with the reduced thin-film energy $E^{red}$ since $h_{ac}(m)$ induces the stray field energy \eqref{eqh2}
as in \eqref{defH-12}. However, \eqref{cond_tan} would be physical relevant for a global minimizer {\it only} if boundary vortices were more expensive than both the N\'eel walls and Bloch line contribution. As explained in the above Claim, this assumption is violated in the regime
\eqref{regime_thin}. Therefore, our energy functional is not adapted for studying global minimizers in the regime \eqref{regime_thin}, but rather for metastable states that satisfy
\eqref{cond_tan}.

Recently, the regime $\textrm{$E^{red}$(Micromagnetic boundary vortex) $\ll \, E^{red}$(N\'eel wall) $\ll
  E^{red}$(Bloch line)}$ was investigated in Ignat \& Kn\"upfer \cite{Ignat_Knuepfer} for thin films of circular cross-section. It is stated that the global minimal configuration for that geometry is given by a $360^\circ-$N\'eel wall
  that concentrates around a radius so that it becomes a vortex (the Landau state of a disk) at the mesoscopic level.

\section{Some preliminaries}

The result stated in Theorem \ref{thm_main} is an extension to the
$S^2-$valued magnetizations of the
following compactness result for $S^1-$valued magnetizations obtained by the authors in \cite{IO} :

\begin{thm}(Ignat \& Otto \cite{IO}) \label{thm1} Let $B^n$ be the unit ball in $\RR^n$, $n=2,3$. For every small $\eta>0$,
let $m'_\eta:B^2\to
S^1$ and $h_\eta=(h'_\eta, h_{3,\eta}):B^3\to \RR^3$ be related by $$ 
\int_{B^3} \left(h'_\eta\cdot \nabla \zeta+h_{3, \eta}\frac{\partial
\zeta}{\partial z}\right)\, dxdz=\int_{B^2} \zeta \nabla\cdot
m'_\eta\, dx,\, \forall \zeta\in C^\infty_c(B^3). $$ Suppose that
\be
\label{bun_stray}
 \int_{B^2}
|\nabla \cdot m'_\eta|^2\, dx+\frac{1}{\eta}\int_{B^3} |h_\eta|^2\, dxdz \leq \frac{C}{\eta |\log \eta|},\ee 
for some fixed constant $C>0$.
Then $\{m'_\eta\}_{\eta\downarrow 0}$ is relatively
compact in $L^1(B^2)$ and any accumulation point $m':B^2\to \RR^2$
satisfies $$|m'|=1 \textrm{ a.e. in $B^2 \quad$ and$ \quad \nabla
\cdot m'=0$ distributionally in $B^2$}.$$
\end{thm}

\bigskip

In the proof of Theorem \ref{detect_vort}, we will use the
following result due to Jerrard \cite{Jer99} for the concentration
of the GL energy \eqref{GLdef} around vortices (see
also Sandier \cite{Sa}, Lin \cite{Lin}):

\begin{thm}(Jerrard \cite{Jer99})
\label{thm_Jer} Let $\alpha\in [0,1)$ and $d>0$ be a
positive integer. There exists $\eps_0=\eps_0(d, \alpha)>0$ such that
for every $0<\eps<\eps_0$, if $m':{\Omega}\to \RR^2$ satisfies
the following conditions:
\be
\label{hypo11}
|m'|\geq \frac 1 2 \quad \textrm{ on } \{x\in \Omega\, :\,
\dist(x, \partial \Omega)\leq r^*(\eps)\} \, \textrm{ for some } r^*(\eps)\in (\frac{1}{|\log \eps|^{4}},1),\ee
$$|\degr(m', \partial \Omega)|=d$$
and
$$\int_\Omega g_\eps(m')\, dx\leq 2\pi (d+\alpha)|\log \eps|,$$ then there exist $n$ points $x_1, \dots, x_n\in \Omega$ with $\dist(x_j, \partial
\Omega)>r^*(\eps)$, $j=1, \dots, n$ and positive integers $d_1,\dots,
d_n>0$ such that the $n$ balls $\{B(x_j, r^*(\eps))\}_{1\leq j\leq n}$ are disjoint,
$$\sum_{j=1}^n d_j=d$$
and
$$\noindent \int_{B(x_j, r^*(\eps))}g_\eps(m')\, dx\geq 2\pi d_j |\log \frac{r^*(\eps)}{\eps}|-C(d,\alpha), \, j=1, \dots, n,$$ where $C(d,
\alpha)$ is a constant only depending on $d$ and $\alpha$.
\end{thm}

In the above theorem, $\Omega$ is any open bounded set (without any regularity condition imposed for the boundary $\partial \Omega$). This is due to hypothesis \eqref{hypo11} of having a security region around $\partial \Omega$.
By degree of a $C^1-$function $v:{\cal C}\to S^1$ defined on a closed curve ${\cal C}$ with the unit tangential vector $\tau$, we mean the winding number
$$\degr(v, {\cal C})=\frac{1}{2\pi} \int_{\cal C} \det(v, \partial_\tau v)\, d\h^1.$$
If $m':{\cal C}\to \RR^2$ is a $C^1-$function with $|m'|>0$ on ${\cal C}$, we set $\degr(m', {\cal C}):=\degr(\frac{m'}{|m'|}, {\cal C})$.
The notion of degree can be extended to continuous vector fields and more generally, $VMO$ vector fields, in particular $H^{1/2}({\cal C}, S^1)$ maps (see Brezis \& Nirenberg \cite{BrNi}).


\section{Proof of Theorem \ref{detect_vort} and Proposition \ref{exec1}}
\label{section5}
First of all, let us define the security region around $\partial \Omega$ together with some notations that we use in the sequel:

\begin{df}
Let $\Omega$ is a simply-connected bounded domain of $C^{1,1}$ boundary. The {\it security region} around $\partial \Omega$
is the maximal set of points around $\partial \Omega$ (in the interior and outside $\Omega$) covered by the normal lines at $\partial \Omega$ before any
crossing occurs. We call {\it depth of the security region} to be the smallest distance to the boundary $\partial \Omega$ where
a crossing of two normal lines occurs and it will be denoted by $R(\partial \Omega)$.
\end{df}

Let $R=R(\partial \Omega)$ be the depth of the security region around $\partial \Omega$.
For $r\in (0, R)$, we denote the interior subdomain $\Omega_r\subset \Omega$ at a distance $r$ from the boundary, i.e.,
\be
\label{def_subdom}
\Omega_r=\{x\in \Omega\,:\, \dist(x, \partial \Omega)> r\} \quad
\textrm{and} \quad \partial \Omega_r=\{x\in \Omega\,:\, \dist(x, \partial \Omega)=r\}
\ee
 be the boundary of this subdomain. For $r\in (-R, 0)$, we write $\partial \Omega_r$ to be the symmetry of $\partial \Omega_{-r}$
across the boundary $\partial \Omega=\partial \Omega_0$ and $\Omega_r\supset \Omega$ be the extended domain surrounded by $\partial \Omega_r$.

Let $l=\h^1(\partial \Omega)$ be the length of $\partial \Omega$.
Set $w:[0,l]\to \partial \Omega$ be a $C^{1,1}$ parametrization of $\partial \Omega$ such that $|\wa|=1$ with
$\wa=\frac{dw}{ds}(s)$ and let $\na=\wa^\perp$ be the outer unit normal vector on $\partial \Omega$ at $w(s)$. Since $\ta=\frac{d^2w}{ds^2}(s)$ is parallel to $\na$ for a.e. $s\in[0, l]$, we will always write $$\ta=\ddot{w}(s) \na$$ where $\ddot{w}(s)$ is the signed length of the vector $\ta$ with respect to $\na$. Notice that $|\ddot{w}(s)|\leq \frac{1}{R(\partial \Omega)}$.
In the security region around $\partial \Omega$, a point $x$ writes in the new coordinates as:
\be
\label{writi}
x=F(s,t)=w(s)+t \na, \quad s\in [0, l], t\in (-R(\partial \Omega), R(\partial \Omega)).\ee 
Note that for interior points $x\in \Omega$, the corresponding normal coordinate $t$ is negative.
We define the symmetry transform $\Phi$ in the security region around $\partial \Omega$:
\be
\label{equst}
\Phi(F(s,t))=F(s,-t) \quad s\in [0, l], t\in (-R(\partial \Omega), R(\partial \Omega)).\ee

\medskip

\bigskip

A first ingredient that we need in the proof of Theorem \ref{detect_vort} is a mirror-reflection extension across the boundary $\partial \Omega$. 
\begin{figure}[htbp]
\center
\includegraphics[scale=0.3,
width=0.3\textwidth]{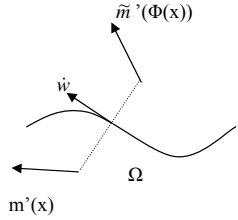} \caption{Mirror-reflection extension.} \label{sym}
\end{figure}
\begin{lem}
\label{lemmir}
Let $R_\infty>0$. There exists $\eps_0=\eps_0(R_\infty)>0$ such that
for every $0<\eps<\eps_0$, the following holds:

Let $\Omega$ be a simply-connected bounded domain of $C^{1,1}$ boundary with the depth of the security region $R(\partial \Omega)\geq R_\infty$. Let $\Phi$ be the symmetry transform across the boundary $\partial \Omega$ defined in \eqref{equst}. In the security region, we consider the interior curve
$$\gamma=\partial \Omega_{\frac{1}{|\log \eps|}}$$
(see notation \eqref{def_subdom})
and
$m':{\Omega}\to \overline{B^2}$ is a $C^1$ vector field that satisfies \eqref{cond_tan},
$$|m'|\geq 1/2 \textrm{ on } \gamma \,\, \textrm{ and }\, \, \deg(m', \gamma)=0.$$ Then
there exists an extension vector field
$\tilde m':\Omega_{-\frac{1}{|\log \eps|}}\to \RR^2$ of $m'$ into the extended domain $\Omega_{-\frac{1}{|\log \eps|}}\supset \Omega$ of boundary $$\tilde \gamma=\Phi(\gamma)=
\partial \Omega_{-\frac{1}{|\log \eps|}}
$$ such that
$$\tilde m'\equiv m' \, \, \textrm{ in } \Omega, \, \,  |\tilde m'|\geq 1/2 \textrm{ on } \tilde \gamma \,\, \textrm{ and }\, \,\deg(\tilde m', \tilde \gamma)=2,$$
\be
\label{in2}
\left|\int_{\tilde \gamma}g_\eps(\tilde m')\, d\h^1(y)-\int_{\gamma}g_\eps(m')\, d\h^1(x)\right|\leq C\bigg({|\log \eps|}\h^1(\partial \Omega)+\frac{1}{|\log \eps|}\int_{\gamma}g_\eps(m')\, d\h^1(x)\bigg)\ee
and
\be
\label{in1}
\left|\int_{\Phi(W)}g_\eps(\tilde m')\, dy-\int_{W}g_\eps(m')\, dx\right|\leq C\bigg(\h^1(\partial \Omega)+\frac{1}{|\log \eps|}\int_{W}g_\eps(m')\, dx\bigg)\ee
where $W\subset \Omega\setminus \Omega_{\frac{1}{|\log \eps|}}$ is any open subset of $\Omega$
and $C=C(R_\infty)$ is a positive constant depending only on $R_\infty$.
\end{lem}

\proof{ of Lemma \ref{lemmir} } We use the notations introduced at the beginning of this section. We have that $|\ddot{w}(s)|\leq \frac{1}{R(\partial \Omega)}\leq \frac{1}{R_\infty}$.
Moreover, differentiating \eqref{writi}, we have that for a.e. $s\in[0, l]$ and $t\in (-R(\partial \Omega), R(\partial \Omega))$,
\be
\label{defF}
DF(s,t)=\left(\begin{array}{cc}\al \wa & \na \end{array} \right) \textrm{ and } DF^{-1}(s,t)=\left(\begin{array}{cc}\frac{1}{\al}\wa\ & \na \end{array} \right)^T,\ee
where $$\al:=1-t\ddot{w}(s).$$
By \eqref{equst} and \eqref{defF}, we compute that:
\be
\label{t1}
S_s(t):=D\Phi(x)=\frac{2}{\al}\wa\otimes \wa-Id \quad \textrm{for a.e. $s\in[0, l]$ and $t\in (-R(\partial \Omega), R(\partial \Omega))$}.\ee
The matrix $S_s(t)$ is symmetric and its inverse is given by $S_s(t)^{-1}=S_s(-t)$.
The mirror-reflection extension $\tm$ of $m'$ is defined as:
\be
\label{t2}
\tmp:=S_s(0)m'(x)=2m'(x)\cdot \wa \wa-m'(x)\,\, \textrm{ for }\, \, x\in \Omega\setminus \Omega_{R(\partial \Omega)}.\ee
(We use that $ a\otimes b \, c=(b\cdot c) a$, for any $a, b, c\in \RR^2$. ) Remark that the condition \eqref{cond_tan} implies that the mirror-reflection extension does not induce jumps at the boundary. Moreover, $|\tmp|=|m'(x)|$ since $S_s(0)=2 \wa\otimes \wa-Id$ is a reflection matrix (i.e., it is symmetric and orthogonal). Therefore, $|\tilde m'|\geq 1/2$ on $\tilde \gamma$.

The goal is to estimate the energies $\ds \int_{\Phi(W)}g_\eps(\tilde m')\, dy$ and $\ds \int_{\tilde \gamma}g_\eps(\tilde m')\, d\h^1$. We start by computing the Dirichlet energy of the extension $\tm$. For that,
we differentiate \eqref{t2} in the coordinates $(s,t)$:
$$D\bigg(\tmp\bigg)
=S_s(0)Dm'(x)DF(s,t)+2\left(\begin{array}{cc} V(s)m'(x) & 0 \end{array} \right),
$$
where
\be
\label{t3}
V(s):=\wa\otimes \ta+\ta\otimes \wa .\ee
Since $D\bigg(\tmp\bigg)=D\tmp D\Phi(x) DF(s,t)$, multiplying by $DF(s,t)^{-1}S_s(-t)$, it implies that
$$D\tmp\stackrel{\eqref{defF}, \eqref{t1}}{=}S_s(0)Dm'(x)S_s(-t)+\frac{2}{\aln}V(s)m'(x)\otimes \wa.$$
Since
$$(D\tmp)^T=S_s(-t)Dm'(x)^T S_s(0)+\frac{2}{\aln}\wa \otimes V(s)m'(x),$$
it follows that
\begin{align}
\nonumber|D\tmp|^2&=tr(D\tmp D\tmp^T)\\
\nonumber&=tr(S_s(0)Dm'(x)S_s(-t)^2Dm'(x)^T S_s(0))+\frac{4}{\aln^2} |V(s)m'(x)|^2\\
\nonumber&\quad+\frac{4}{\aln} tr(S_s(0)Dm'(x)S_s(-t)\wa \otimes V(s)m'(x))\\
\label{descomp}&=I+II+III.
\end{align}
For the first term in \eqref{descomp}, we compute that
$$S_s(-t)^2\stackrel{\eqref{t1}}{=}-\frac{4t\ddot{w}(s)}{\aln^2} \wa \otimes \wa+Id.$$
Since $tr(SAS^{-1})=tr(A)$ and $tr(Av\otimes Av)=|Av|^2\leq |A|^2 |v|^2$ for any two matrices $A$ and $S$ in $\RR^{2\times 2}$ with $S$ invertible and any vector $v\in \RR^2$, we deduce that
\be
\label{e1}
I\leq \left(1+\frac{4|t||\ddot{w}(s)|}{\aln^2} \right)|Dm'(x)|^2.\ee
For the second term in \eqref{descomp}, we have that
$|V(s)|^2\stackrel{\eqref{t3}}{=}2|\wa \otimes \ta|^2=2|\ddot{w}(s)|^2$ and therefore,
\be
\label{e2}
II\leq \frac{4}{\aln^2}|V(s)|^2|m'(x)|^2 \leq \frac{8|\ddot{w}(s)|^2}{\aln^2}|m'(x)|^2.\ee
For the third term in \eqref{descomp}, we compute that
$$S_s(-t)\wa\stackrel{\eqref{t1}}{=}\frac{\al}{\aln} \wa \quad \textrm{ and }\quad S_s(0)V(s)\stackrel{\eqref{t1}}{=}\ddot{w}(s)\left( \begin{array}{cc} 0 & 1\\ -1 & 0 \end{array} \right).$$
Using that $tr(A b\otimes c)=c\cdot Ab$ and $tr(A)=tr(S_s(0)AS_s(0))$ for any matrix $A$ in $\RR^{2\times 2}$ and any vectors $b,c\in \RR^2$, we deduce that
\begin{align}
\nonumber III&=\frac{4\al}{\aln^2} tr\bigg(Dm'(x)\wa \otimes (S_s(0)V(s)m'(x))\bigg)\\
\nonumber&=-\frac{4\al \ddot{w}(s)}{\aln^2} tr(Dm'(x)\wa \otimes m'(x)^\perp)\\
\nonumber&=-\frac{4\al \ddot{w}(s)}{\aln^2} m'(x)^\perp\cdot (Dm'(x)\wa)\\
\label{e3}&\leq \frac{4\al |\ddot{w}(s)|}{\aln^2} |m'(x)|\, |Dm'(x)|.
\end{align}
Since $|\det D\Phi(x)|\stackrel{\eqref{t1}}{=}\frac{\aln}{\al}$, we deduce by \eqref{descomp}, \eqref{e1}, \eqref{e2} and \eqref{e3},
\begin{align}
\nonumber |D\tmp|^2 |\det(D\Phi(x))|&\leq \frac{\aln}{\al} \left(1+\frac{4|t||\ddot{w}(s)|}{\aln^2} \right)|Dm'(x)|^2+\frac{8|\ddot{w}(s)|^2}{\aln \al}|m'(x)|^2\\
\label{calfin} & \quad + \frac{4|\ddot{w}(s)|}{\aln} |m'(x)|\, |Dm'(x)|.\end{align}
Therefore, for every open set $W\subset  \Omega\setminus \Omega_{\frac{1}{|\log \eps|}}$, we obtain by Young's inequality,
\begin{align}
\nonumber
\int_{\Phi(W)} |D\tm(y)|^2\, dy&=\int_W |\det D\Phi(x)|\, |D\tmp|^2\, dx\\
\nonumber
&\stackrel{\eqref{calfin}}{\leq}\int_W \bigg\{ (1+\frac{C}{|\log \eps|}) |Dm'(x)|^2+ C{|\log \eps|}|m'(x)|^2\bigg\}\, dx\\
\nonumber
&\leq(1+\frac{C}{|\log \eps|}) \int_W |Dm'(x)|^2\, dx +C \h^1(\partial \Omega),
\end{align}
with $C=C(R_\infty)>0$ and $\eps\leq \eps(R_\infty)$.
(We use that $\h^2(W)\leq \frac{C}{|\log \eps|} \h^1(\partial \Omega)$.)
Also,
\begin{align*}
\int_{\Phi(W)} (1-|\tm(y)|^2)^2\, dy&=\int_W |\det D\Phi(x)|\,(1-|m'(x)|^2)^2 \, dx\\
&\leq (1+\frac{C}{|\log \eps|}) \int_W (1-|m'(x)|^2)^2 \, dx.\end{align*}
Therefore, we obtain:
$$\int_{\Phi(W)} g_\eps(\tm)\, dy\leq(1+\frac{C}{|\log \eps|})\int_{W} g_\eps(m')\, dx +C \h^1(\partial \Omega).$$
By changing $t$ to $-t$ in the above argument, the inverse inequality also holds:
\be
\label{bunineq1}
\int_{W} g_\eps(m')\, dx\leq (1+\frac{C}{|\log \eps|}) \int_{\Phi(W)} g_\eps(\tm)\, dy +C \h^1(\partial \Omega).
\ee
Thus, inequality \eqref{in1} immediately follows.
For proving inequality \eqref{in2}, we proceed in the same way: Since
$F(\cdot, -\frac{1}{|\log \eps|})=w-\frac{1}{|\log \eps|}{\bf \nu}$ is a Lipschitz parametrization of $\gamma$, we compute
$$\bigg|\frac{d}{ds}(\Phi(F(s, -\frac{1}{|\log \eps|})))\bigg|\stackrel{\eqref{t1}}{=}|\det D\Phi(F(s, -\frac{1}{|\log \eps|}))|\bigg|\frac{d}{ds}(F(s, -\frac{1}{|\log \eps|}))\bigg|$$ and we have by \eqref{calfin},
\begin{align}
\nonumber
\int_{\Phi(\gamma)} g_\eps(\tm(y))\, d\h^1(y)&=\int_{\gamma} |\det D\Phi(x)|\, g_\eps(\tm(\Phi(x)))\, d\h^1(x)\\
\nonumber
&\leq(1+\frac{C}{|\log \eps|}) \int_\gamma \, g_\eps(m'(x))d\h^1(x) +C{|\log \eps|} \h^1(\partial \Omega).
\end{align}
By symmetry, \eqref{in2} follows immediately.

It remains to prove that if $\deg(m', \gamma)=0$, then
$\deg(\tilde m', \tilde \gamma)=2$. For that let $\f_0:[0,l]\to \RR$ be the lifting of ${\bf \dot{w}}$, i.e., $\wa=e^{i\f_0(s)}$. Obviously, \be
\label{scriu2} \deg({\bf \dot{w}})=
\frac{1}{2\pi}(\f_0(l)-\f_0(0))=1.\ee 
On the curve $\gamma$, we know that $m'\in C^1(\gamma, \RR^2)$ and we write $m'=\rho v$ with $\rho=|m'|\geq 1/2$ and $v:\gamma\to S^1$. Then $\rho, v\in C^1(\gamma)$. Then $\deg(v, \gamma)=\deg(m', \gamma)=0$. In this case, the theory of lifting yields the existence of a lifting $\f\in C^1(\gamma, \RR)$ such that
$v=e^{i\f}$. If $t:=\frac{1}{|\log \eps|}$, then $F(\cdot, -t)$ is a parametrization of $\gamma$ and we have
\be
\label{scriu_lift}
0=\deg(v, \gamma)=\frac{1}{2\pi}(\f(F(l,-t))-\f(F(0,-t))).
\ee
Notice that the reflection matrix $S_s(0)$ has the following form: $$S_s(0)\stackrel{\eqref{t1}}{=}\left(\begin{array}{cc} \cos 2\f_0(s) & \sin 2\f_0(s)\\  \sin 2\f_0(s) & -\cos 2\f_0(s) \end{array}\right).$$ That implies the following writing of $\tilde m'$ on the curve $\tilde \gamma=\Phi(\gamma)$ parametrized by $F(\cdot, t)$:
$$\tilde m'(F(s,t))=\rho(F(s,-t))S_s(0)v(F(s,-t))=\rho(F(s,-t)) e^{i\bigg(2\f_0(s)-\f(F(s,-t)\bigg)}.$$
Therefore, by \eqref{scriu2} and \eqref{scriu_lift}, we conclude that
$$\deg(\tilde m', \tilde \gamma)=2.$$
\qed

We now prove the concentration of the Ginzburg-Landau energy on a small region (either one interior vortex, or two boundary vortices) under
the condition \eqref{cond_tan} in the regime \eqref{en_vort}:

\proof{ of Theorem \ref{detect_vort}} Let $R=R(\partial \Omega)$ be the depth of the security region around $\partial \Omega$.
We proceed in several steps:

\bigskip

\noindent{\bf Step 1.} {\it Find a good set of boundaries.} We define
the set $I$ of distances $r\in(\eps, R)$ such that we control the energy of $m'$ on the
boundary $\partial \Omega_r$ (and consequently, the modulus $|m'|$ via \eqref{modmic}), i.e., \be \label{circles} I=\left\{r\in
(0, R-\eps)\, :\, \int_{\partial \Omega_r} g_\eps(m')\,
d\h^1\leq {|\log \eps|^3}\right\}.\ee How large is
the set $I$? We show that for each interval $J\subset (0, R-\eps)$ of
length $\ell \gg \frac{1}{|\log \eps|^2}$, there exist infinitely many distances $r$ belonging
to $I\cap J$. More precisely, we have for small $\eps>0$ that
\be
\label{eqre}
\h^1(I\cap J)\geq \frac{\ell}{2}.\ee Indeed,
one has
$$4\pi |\log \eps|\stackrel{\eqref{en_vort}}{\geq} \int_{J \setminus I}\int_{\partial \Omega_r} g_\eps(m')\,d\h^1
dr\geq |\log \eps|^3 \h^1(J \setminus I)
$$ which yields $\h^1(J\setminus I)\leq \frac{4\pi}{|\log \eps|^2}\leq \frac{\ell}{2}$ for small $\eps>0$
and therefore, \eqref{eqre} holds.
Moreover, remark that 
$
|m'|\geq \frac 1 2$ for every $r\in I$,
if $\eps>0$ is small enough.
Indeed, since $r< R-\eps$ it means that $\h^1(\partial \Omega_r)\geq \h^1(\partial B(0,\eps))\geq \eps$. Denoting by $\rho:=|m'|$ and $min:=\min \{\rho(x)\, :\, x \in \partial \Omega_r\}$, it is easy to check that
(see Lemma 2.3. in \cite{Jer99}):
\be
\label{modmic}
{|\log \eps|^3}\geq  \int_{\partial \Omega_r} g_\eps(m')\,
d\h^1\geq \int_{\partial \Omega_r} \bigg(|\partial_\tau \rho|^2+ \frac{1}{\eps^2} (1-\rho^2)^2 \bigg)\, d\h^1\geq \frac{C}{\eps}(1-min)^2,\ee
where $\tau$ is the tangent unit vector at $\partial \Omega_r$. Thus, one concludes that $(1-min)^2\lesssim \eps |\log \eps|^3\ll1$, i.e., $min\geq 1/2$ for small $\eps>0$.
(Relation \eqref{modmic} is obvious if $\rho$ is constant (equal with $min$). Otherwise, the GL energy of the modulus $\rho$ controls the following quantity $\frac{1}{\eps}\int_{Im(\rho)}(1-
y^2)\, dy$ on the image set $Im(\rho)$ of $\rho$ and forces $\rho$ to take values close to $1$.)

\bigskip

By\eqref{eqre}, we can choose $r_1\in I\cap (\frac{1}{2|\log \eps|}, \frac{2}{|\log \eps|} )$. W.l.o.g., we may suppose that $$r_1=\frac{1}{|\log \eps|}.$$
We distinguish two cases in function of $\degr(m',
\partial \Omega_{r_1})$.
\bigskip

\noindent {\bf Step 2}. {\it We assume that $|\degr(m',
\partial \Omega_{r_1})|> 0.$} W.l.o.g. we may suppose that $d:=\degr(m',
\partial \Omega_{r_1})\geq 1$.

\bigskip

\noindent {\it Extension.} We will extend the vector field $m'\bigg|_{\Omega_{r_1}}$ by a vector field $\tilde m'$ defined on the larger domain $\Omega_{r_1-r^*}\supset \Omega_{r_1}$ with
$$r^*=\frac{1}{|\log \eps|^3}\in (0, \frac{r_1}{2})$$ such that $\tilde m'=m'$ in $\Omega_{r_1}$ and
\begin{align}
\nonumber &\degr(\tilde m',
\partial \Omega_{r_1-r^*})=d\geq 1,\\
\nonumber &|\tilde m'|\geq \frac 1 2 \quad \textrm{ in }
\Omega_{r_1-r^*} \setminus \Omega_{r_1}\\
&\label{buna_en_lev} \int_{\Omega_{r_1-r^*}} g_{\eps}(\tilde m')\, dx\leq
2\pi(1+\alpha)|\log \eps|+C.\end{align}
For that, using the notation \eqref{writi}, for each point $x=F(s,t)\in \Omega_{r_1-r^*} \setminus \Omega_{r_1}$ (here, $t<0$), we consider $y_x=F(s,-r_1)\in \partial
\Omega_{r_1}$ to be the normal projection of $x$ on $\partial
\Omega_{r_1}$ with $\dist(x,
\partial \Omega_{r_1})=|x-y_x|=|t+r_1|$. Then we define $\tilde m':\Omega_{r_1-r^*} \to {\RR^2}$ as
$\tilde m'=m'$ in $\Omega_{r_1}$ and \be \label{defmtild} \tilde
m'(x):=m'(y_x),\quad \textrm{ for every } x\in \Omega_{r_1-r^*} \setminus
\Omega_{r_1}.\ee Since $r_1\in I$, \eqref{modmic} implies $|\tilde m'|\geq 1/2$ on
$\Omega_{r_1-r^*} \setminus \Omega_{r_1}$, $\deg(\tilde m',
\partial \Omega_{r_1-r^*})=d\geq 1$ and
\begin{align}
\nonumber \int_{\Omega_{r_1-r^*} \setminus \Omega_{r_1}} g_{\eps}(\tilde m')\,
dx&=\int_{r_1-r^*}^{r_1}\int_{\partial \Omega_r} g_{\eps}(\tilde m') \, d\h^1 dr\\
\nonumber &\stackrel{\eqref{defmtild}}{\leq} C r^*
\int_{\partial \Omega_{r_1}}g_{\eps}(m') \, d\h^1\\
\label{con24}&\stackrel{\eqref{circles}}{\leq}
C r^* |\log \eps|^3= C.\end{align}
Thus, by \eqref{en_vort}, we
obtain \eqref{buna_en_lev}. By Theorem \ref{thm_Jer} and \eqref{buna_en_lev}, we deduce that $d=1$ and there exists
a point $x_1\in \Omega_{r_1}$ such that
$$\noindent \int_{B(x_1, r^*)}g_\eps(\tilde m')\, dx
\geq 2\pi |\log \frac{r^*}{\eps}|-C(\alpha),$$
where $C(\alpha)$ is a constant depending only on
$\alpha$. Therefore, we obtain via \eqref{con24} that
$$ \int_{B(x_1, r^*)}g_{\eps}(m')\, dx\geq
2\pi |\log\frac{r^*}{\eps}|-C(\alpha).$$ Since $B(x_1, r^*)\subset \Omega$, the conclusion
\eqref{estvoin} follows.

\bigskip

\noindent {\bf Step 3}. {\it We now deal with the other case $\degr(m',
\partial \Omega_{r_1})=0.$}

\bigskip

\noindent {\it Mirror-reflection extension.} We consider the
symmetry transform $\Phi$ defined in Lemma \ref{lemmir} across the
boundary $\partial \Omega$ together with the mirror-reflection extension
$\tilde m': \Omega_{r_2}\to \overline{B^2}$ where $r_2=-{r_1}$. Then Lemma
\ref{lemmir} yields \be \label{con78}
\int_{\Omega_{r_2}} g_{\eps}(\tilde m')\, dx\leq 2\pi
(2+2\alpha)|\log \eps|+C(\partial \Omega),\ee $|\tilde m'|\geq \frac 1 2$ on
$\partial \Omega_{r_2}$, the
degree of $\tilde m'$ on the boundary $\partial \Omega_{r_2}$ is
equal to $2$ and
$$
\int_{\partial \Omega_{r_2}} g_{\eps}(\tilde m')\, d\h^1\leq
2|\log \eps|^3$$
for $\eps$ small enough.
The extension argument in Step 2 leads via Theorem \ref{thm_Jer} to the concentration of
the Ginzburg-Landau energy of $\tilde m'$ into two vortex balls $B(x_2, r^*)$
and $B(x_3, r^*)$ with $x_2, x_3\in \Omega_{r_2}$ and there exist two non-negative numbers $d_2\geq d_3\geq 0$, $d_2+d_3=2$ such that
\be
\label{es12}
\int_{B(x_j, r^*)}g_{\eps}(\tilde m')\, dx\geq
2\pi d_j |\log\frac{r^*}{\eps}|-C,\quad  j=2,3.\ee
(The assumption
$\alpha<1/2$ is needed so that $2+2\alpha<3$.)

{\it Case 1}: $d_2=2$ (i.e., there is one vortex ball of degree $2$ in $\Omega_{r_2}$). The level of energy \eqref{en_vort} rules out that $B(x_2, r^*)\subset \Omega$. By Lemma \ref{lemmir}, it also means that
$B(x_2, r^*)$ is not included in $\Omega_{r_2}\setminus \Omega$ (otherwise, the symmetry of the energy distribution around the boundary would imply again that the reflected domain
$\Phi(B(x_2, r^*))$ charges the energy more than the level \eqref{en_vort} in the interior of $\Omega$). Therefore, $B(x_2, r^*)\cap \partial \Omega\neq \emptyset$. Choose $x^*_2=x^*_3\in B(x_2, r^*)\cap \partial \Omega$. Then $B(x_2, r^*), \Phi(B(x_2, r^*))\subset B(x^*_2, 10r^*)$ and by Lemma \ref{lemmir} and \eqref{es12}, we conclude that
\begin{align*}
\int_{B(x^*_2, 10r^*)\cap \Omega}g_{\eps}(m')\, dx&\geq \frac 1 2 \left(\int_{B(x_2, r^*)\cap \Omega} g_{\eps}(\tilde m')\, dx+\int_{\Phi(B(x_2, r^*)\setminus \Omega)} g_{\eps}(\tilde m')\, dx\right)\\
&\stackrel{Lemma\ref{lemmir}}{\geq} \frac 1 2 \left(\int_{B(x_2, r^*)\cap \Omega} g_{\eps}(\tilde m')\, dx+\int_{B(x_2, r^*)\setminus \Omega} g_{\eps}(\tilde m')\, dx\right)-C\\
&{\geq} \frac 1 2 \int_{B(x_2, r^*)}g_{\eps}(\tilde m')\, dx-C\stackrel{\eqref{es12}}{\geq} 2\pi |\log\frac{r^*}{\eps}|-C.
\end{align*} Here, \eqref{estvobd} holds and $B(x^*_2, 10r^*)=B(x^*_3, 10r^*)$ are the boundary vortex balls.

{\it Case 2}:  $d_2=d_3=1$ (i.e., there are two disjoint vortex balls of degree $1$ in $\Omega_{r_2}$). If \eqref{estvoin} holds, then we are done. Suppose that \eqref{estvoin} is not satisfied. Then we want to prove
\eqref{estvobd}.
As in Case 1, the symmetry of the energy distribution around the boundary implies via Lemma \ref{lemmir} that
none of the balls $B(x_2, r^*)$
and $B(x_3, r^*)$ is included in $\Omega$ or $\Omega_{r_2}\setminus \Omega$ (otherwise, \eqref{estvoin} would hold). Therefore,
$B(x_j, r^*)\cap \partial \Omega\neq \emptyset$ for $j=2,3$.
Choose $x^*_j\in B(x_j, r^*)\cap \partial \Omega$ for $j=2,3$. Then
$B(x_j, r^*), \Phi(B(x_j, r^*))\subset B(x^*_j, 10r^*)$ for $j=2,3$.
As before, by Lemma \ref{lemmir} and \eqref{es12}, we conclude that
\begin{align*}
\int_{(B(x^*_2, 10r^*)\cup B(x^*_3, 10r^*))\cap \Omega}g_{\eps}(m')\, dx
& \geq \frac 1 2 \bigg(\int_{\big(B(x_2, r^*)\cup B(x_3, r^*)\big) \cap \Omega} g_{\eps}(\tilde m')+\\&\quad \quad \quad \quad \int_{\Phi\big(B(x_2, r^*)\setminus \Omega\big)\cup 
\Phi\big(B(x_3, r^*)\setminus \Omega\big)} g_{\eps}(\tilde m')\bigg)\\
&\geq\frac 1 2 \int_{B(x_2, r^*)\cup B(x_3, r^*)}g_{\eps}(\tilde m')\, dx-C\\&\geq  2\pi |\log\frac{r^*}{\eps}|-C.
\end{align*}
(Here, we used that $\Phi\bigg(B(x_2, r^*)\setminus \Omega\bigg)\cap \Phi\bigg(B(x_3, r^*)\setminus \Omega\bigg)=\emptyset$ since the two balls $B(x_2, r^*)$
and $B(x_3, r^*)$ are disjoint and lie in the security region of $\partial \Omega$.)
\qed

\bigskip

The natural question is whether the lower bound for the energy of a boundary vortex given in Theorem \ref{detect_vort} is optimal. A positive answer is supported by the following result:
we prove that
the loss of energy of order $\frac{|\log \eps|}{\log |\log \eps|}$ (with respect to $\pi |\log \eps|$ which is the exact half energy of an interior vortex) may be achieved for $C^1$ domains.

\proof{ of Proposition \ref{exec1}} 
The aim is to construct a boundary vortex on $\partial \Omega$ centered at the origin. 

\noindent {\bf Step 1}. {\it Construction of $m'_\eps$}. We define the two-dimensional vector field $m'_\eps$ that is tangent at $\partial \Omega$ and its phase $\ph$ is linear on every arc of circle $\{|x|=r\}\cap \Omega$ with $r \in (0, 1/200)$. More precisely, let
$$m'_\eps(x)=\begin{cases}
e^{i\ph(x)} & \quad \textrm{if } x\in \Omega  \textrm{ and  } \eps<|x|<1/200,\\
\frac{|x|}{\eps}e^{i\ph(x)} & \quad \textrm{if } x\in \Omega  \textrm{ and  } 0<|x|<\eps,
\end{cases}$$
where the phase $\ph$ is given in the polar coordinates as follows: $\ph(r, \cdot): (-\ga(r), \ga(r))\to (-\pi/2, \pi/2)$ is an odd function (i.e., $\ph(r, \theta)=-\ph(r, -\theta)$) and it is linear in $\theta$,
$$\ph(r, \theta)=\bigg( 1+\frac{\dte(r)}{\ga(r)}\bigg) \theta \quad \textrm{ for every} \quad \theta\in (-\ga(r), \ga(r)), r\in (0, 1/200).$$
The phase correction $\dte:(0, 1/200)\to (-\frac \pi 2,0)$ due to the condition \eqref{cond_tan} is defined as
$$e^{i\dte(r)}=\frac{1+ir \ga'(r)}{\sqrt{1+(r \ga'(r))^2}}  \quad \textrm{ for every} \quad r\in (0, 1/200).$$
(Indeed, one can easily check that $m'_\eps$ is tangent at $\partial \Omega$.)

\noindent {\bf Step 2}. {\it Estimate of the GL energy outside the core region}. We first estimate the energy of $m'_\eps$ away
from the core, i.e., ${\cal D}_1=\{x\in \Omega\, : \,
\eps<|x|<1/200\}$. 
For that, we need the following computations: for $r\in (0, 1/200)$,
\be
\label{com1}
\ga'(r)=\frac{-1}{r \log \frac 1 r (\log \log \frac 1 r)^2},
\ee
\be
\label{com2}
\dte(r)=-\arccos \frac{1}{\sqrt{1+(r \ga'(r))^2}} \Rightarrow |\dte(r)|\leq 2|r\ga'(r)|\stackrel{\eqref{com1}}{\leq} \frac{2}{ \log \frac 1 r (\log \log \frac 1 r)^2},
\ee
where we used that $\arccos:[-1,1]\to [0, \pi]$ satisfies $\sqrt{1-t^2}\leq \arccos t \leq 2\sqrt{1-t^2}$ for $t\in [\frac 1 2,1]$.
By a change of variable $r=e^s$ and differentiating \eqref{com2} in the new logarithmic variable $s$, we deduce that
\be
\label{line1}
r\bigg|\frac{d}{dr}\dte(r)\bigg|=\bigg|\frac{d}{ds}\dte(e^s)\bigg|=\frac{1}{1+\bigg(\frac{d}{ds}\ga(e^s)\bigg)^2} \bigg|\frac{d^2}{ds^2}\ga(e^s)\bigg|
\ee and
we have for a universal constant $C>0$, 
 \begin{align*} \int_{{\cal
D}_1}g_{\eps}(m'_\eps)\, dx&=\int_{{\cal D}_1} |\nabla \ph|^2\, dx&\\
&=
\int_\eps^{1/200}\int_{-\ga(r)}^{\ga(r)}\bigg(
\frac{|\partial_\theta \ph|^2}{r}+r{|\partial_r \ph|^2}\bigg)\,
d\theta dr\\
&\stackrel{\eqref{com1}, \eqref{com2}}{\leq} \int_\eps^{1/200}
\bigg\{\frac{2\ga(r)}{r} \bigg(1+\frac{C}{ \log \frac 1 r (\log \log \frac 1 r)^2}\bigg)
+C r \ga(r) \bigg(\frac{1}{ (\log r)^2 r^2}+ (\frac{d}{dr}\dte(r))^2\bigg)\bigg\}\, dr\\
&{\leq} \int_\eps^{1/200}
\bigg\{2\ga(r)+C \ga(r) r^2 (\frac{d}{dr}\dte(r))^2\bigg\}\, \frac{dr}{r}+C\\
&\stackrel{s=\log r, \, \eqref{line1}}{\leq} \int_{\log \eps}^{\log \frac{1}{200}} \bigg\{
2\ga(e^s)+\ga(e^s) 
\frac{C}{\bigg(1+(\frac{d}{ds}\ga(e^s))^2\bigg)^2} \bigg(\frac{d^2}{ds^2}\ga(e^s)\bigg)^2 \bigg\}
\,
ds+C\\
&\leq \int_{\log \eps}^{\log \frac{1}{200}} \bigg( \pi 
-\frac{2}{\log |s|}\bigg)\bigg(1+C\bigg(\frac{d^2}{ds^2}\ga(e^s)\bigg)^2 \bigg)\,
ds+C\\
&\leq \pi |\log \eps|- 2\int_{\log \eps}^{\log \frac{1}{200}} \frac{1}{\log |s|} \bigg(1-\frac{C}{s^4} \bigg)\,
ds+C\\
&\leq (\pi -\frac{C}{\log|\log \eps|}) |\log \eps|.
\end{align*} 
(Here we used that $\int_{10}^x \frac{1}{\log s}\, ds\sim \frac{x}{\log x}$ as $x\to \infty$.)

\noindent {\bf Step 3}. {\it Estimate of the GL energy inside the core region}.
Now we estimate the energy of $m'_\eps$ on the
core, i.e., ${\cal D}_2=\{x\in \Omega\, : \, 0<|x|<\eps\}$. Using the same argument as above and the change of coordinates $r=e^s$,
we compute
\begin{align*}
\int_{{\cal D}_2} |\nabla m'_\eps|^2 \, dx&=
\int_0^{\eps}\int_{-\ga(r)}^{\ga(r)} \bigg(
\frac{r^2|\nabla \ph|^2}{\eps^2}+\frac{1}{\eps^2}\bigg)\, rd\theta dr\\
&\stackrel{\eqref{com1}, \eqref{com2}}{\leq} \int_0^\eps \bigg(
\frac{Cr}{\eps^2}+\frac{r^3\ga(r)}{\eps^2} \bigg(\frac{d}{d r}\dte(r)\bigg)^2\bigg)\, dr\stackrel{\eqref{line1}}{\leq}\int_0^\eps
\frac{Cr}{\eps^2}\, dr =O(1)
\end{align*}  and $$\int_{{\cal D}_2}
\frac{1}{\eps^2}(1-|m'_\eps|^2)^2 \, dx=O(1).$$ 

\noindent {\bf Step 4}. {\it Conclusion}.
The $H^1$ vector field $m'_\eps$ satisfies the required properties in Proposition \ref{exec1}. However, $m'_\eps$ is not $C^1$. In order to construct a smooth 
$m'_\eps$, we define $f:\RR\to [0,1]$ a smooth function such that $f(t)=0$ if $t\leq 0$ and $f(t)=1$ if $t\geq 1$. Now, it is enough to change the vector field $m'_\eps$ defined at Step 1
only in the core region as follows: $m'_\eps(x)=f(\frac{|x|}{\eps})e^{i\varphi_\eps(x)}$ if $0\leq |x|\leq \eps$.
\qed

\begin{rem}
\label{remu} The GL energy of a boundary vortex placed in a corner is
proportional with the corner angle. Therefore, the loss of energy of order $\frac{|\log \eps|}{\log |\log \eps|}$ for $C^1$ boundaries (see Proposition \ref{exec1}) increases to an order of $|\log \eps|$ for Lipschitz boundaries. More precisely, let
$\Omega=\{(x_1, x_2)\,:\, x_1\in(-1, 1), \, x_2> |x_1| \tan
\frac{\pi-\alpha }{2}\}$ where $\alpha\in (0,\pi)$ is the corner
angle of $\Omega$ at the origin. For every $0<\eps<1$, we consider the following
approximation of a vortex:
$$m'_\eps(x)=\begin{cases}
\frac{x}{|x|} & \quad \textrm{if } \eps<|x|<1,\\
\frac{x}{\eps} & \quad \textrm{if } 0<|x|<\eps.\\
\end{cases}$$
Then $m'_\eps$ satisfies \eqref{cond_tan} on $\partial \Omega\cap
B^2$ and
$$\int_{\Omega\cap B^2} g_\eps(m'_\eps)\, dx\leq \alpha |\log \eps|+O(1).$$
\end{rem}


\section{Proof of Theorem \ref{thm_main}}
\label{section6}

We will work at the level of sequences of parameters $\eps_k$ and $\eta_k$ and a sequence of magnetizations $m_k$ satisfying the assumptions in Theorem \ref{thm_main}.
We will prove the Theorem in a slightly larger regime than \eqref{regim_nou}; more precisely, it is enough to assume that 
\be
\label{new_regi}
\eps_k^\beta\lesssim \eta,\ee
for some constant $\beta\in (0,1)$ such that
$$\beta<1-\alpha.$$
By \eqref{en_level2}, let $A>0$ be such that
\be \label{regim_cond} 
E_{\eps_k, \eta_k}(m_k)-2\pi |\log \eps_k|
\leq
\frac{A}{\eta_k |\log
\eta_k|} \quad \textrm{ for every }k\in
\NN.\ee
By \eqref{regim_cond11}, there exists $a>0$ such that 
\be
\label{ect}
a\log |\log \eps_k|\leq \frac{A}{\eta_k |\log
\eta_k|},\ee for every $k\in \NN$.
Let $U_k:\RR^3\to \RR$ be the stray field potential associated to
$m_k$ defined by \eqref{1c} for $(\nabla,\frac{\partial}{\partial
z}) U_k$ that satisfies
\be
\label{choiceUk}
\int_{\RR^2\times \RR}\left(|\nabla U_k|^2+\big|\frac{\partial U_k}{\partial
z}\big|^2 \right)\,dx dz=\frac 1 2 \int_{\RR^2}
\left|\,|\nabla|^{-1/2}(\nabla\cdot m'_k)\right|^2\, dx.\ee
By the Lax-Milgram theorem, the potential $U_k$ exists and is unique in the Beppo-Levi space (see Dautray and Lions
\cite{DauLio}):
$${\cal BL}=\{U:\RR^3\to \RR\,:\, (\nabla,\frac{\partial}{\partial
z}) U\in L^2(\RR^3), \frac{U}{1+|x|}\in L^2(\RR^3)\}.$$

We proceed in several steps:

\noindent {\bf Step 1}. {\it Location of the vortex balls of $m'_k$.} Let $r^*=1/|\log \eps_k|^3$.
By Theorem \ref{detect_vort}, there exist at most two points $x_k, \tilde x_k \in \bar \Omega$
such that
\be \label{con26}
\int_{(B(x_k, r^*)\cup B(\tilde x_k, r^*))\cap \Omega}g_{\eps_k}(m'_k)\, dx\geq 2\pi |\log
\frac{r^*}{\eps_k}|-C\geq 2\pi
|\log{\eps_k}|-100\log|\log \eps_k|,\ee for $k$ sufficiently large.
Obviously, up to a
subsequence, $\{x_k\}, \{\tilde x_k\}\subset \overline \Omega$ converge to two points
$x_0, \tilde x_0\in \overline \Omega$ and we have that for every small $\sigma>0$,
$$B(x_k, r^*)\subset B(x_0, \sigma), \, B(\tilde x_k, r^*)\subset B(\tilde x_0, \sigma)$$ for $k$ sufficiently large.

The set ${\cal D}=B(x_0, \sigma)\cup B(\tilde x_0, \sigma)$ is the location of the essential topological defects of each
$m'_k$. Now the goal is to prove that $\{m_k\}$ is relatively
compact in $L^1(\Omega\setminus {\cal D})$. The idea is to
approximate $m'_k$ away from ${\cal D}$ by $S^1-$valued vector
fields, denoted by $M'_k$ that satisfy the hypothesis of Theorem~\ref{thm1}. For that, let ${\bbbb}\subset
\Omega\setminus {\cal D}$ be an arbitrary ball. To
simplify the notation, let ${\bbbb}=B(0,2)\subset \RR^2$ be
the ball of radius $2$ centered in the origin. Since $m_{3,k}^2\geq m_{3,k}^4=(1-|m'_k|^2)^2$, the energy level on $\bbbb$ is bounded as
follows:
\begin{align} \nonumber
\int_{\bbbb} |\nabla m_k|^2\, dx+\frac{1}{\eps_k^2}\int_{\bbbb}
(1-|m'_k|^2)^2\, dx+&\frac 1 \eta_k \int_{\RR^2}
||\nabla|^{-1/2}(\nabla \cdot m'_k)|^2\,
dx\\
\nonumber &\leq E_{\eps_k, \eta_k}(m_k)-\int_{{\cal D}\cap
\Omega}g_{\eps_k}(m'_k)\, dx\\
\nonumber & \stackrel{\eqref{en_level}, \eqref{regim_cond}, \eqref{con26}}{\leq}
\min \{\frac{A}{\eta_k |\log \eta_k|},  2\pi \alpha |\log \eps_k|\}+100\log |\log \eps_k|\\
\label{en_D} &\stackrel{\eqref{ect}}{\leq} \min\bigg\{ \frac{\tilde{A}}{\eta_k |\log \eta_k|}, 2\pi \alpha |\log \eps_k|+100\log |\log \eps_k|\bigg\}\end{align} for $k$ sufficiently large and $\tilde A=A(a+100)/a$ (by \eqref{ect}).

\medskip

\noindent {\bf Step 2}. {\it Construction of a square grid.} For each shift
$t\in [0, \eps_k^{\beta})$, write
\begin{figure}[htbp]
\center
\includegraphics[scale=0.2,
width=0.2\textwidth]{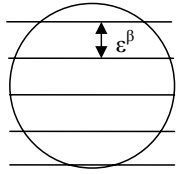} \caption{The net of horizontal lines.} \label{Netul}
\end{figure}
$$V_t:=\{(x_1,x_2)\in {\bbbb}\, :\, x_2\equiv t \quad (\mod
\eps_k^{\beta})\}$$ for the net of horizontal lines at a distance
$\eps_k^{\beta}$ in $\bbbb$. By the mean value theorem, 
there exists $t_k\in ({0,\eps_k^\beta})$
such that $$\int_{V_{t_k}} g_{\eps_k}(m'_k)\, d\h^1\leq \frac{1}{\eps_k^\beta} \int_{\bbbb} g_{\eps_k}(m'_k)\,
dx.$$ 
If one repeats the above argument for the net of vertical
lines at a distance $\eps_k^{\beta}$ in $\bbbb$, we get a square
grid ${\cal R}_k$ of size $\eps_k^\beta$ such that the convex hull of
${\cal R}_k$ covers the unit ball $B^2\subset \bbbb$ and \be
\label{condR}
 \int_{{\cal R}_{k}} g_{\eps_k}(m'_k)\,
d\h^1 \stackrel{\eqref{en_D}}{\leq} \min \{\frac{2\tilde A}{\eps_k^{\beta}\eta_k |\log
\eta_k|}, \frac{C
|\log \eps_k|}{\eps_k^{\beta}}\}.\ee
By the same argument as in \eqref{modmic}, the estimate \eqref{condR} together with $\beta<1$ implies that 
${\cal R}_{k}\subset \{|m'_k|> 1/2 \}$ for $k$ large enough.

\noindent {\bf Step 3}. {\it Vanishing degree on the cells of the grid.}
In order to approximate $m'_k$ in $B^2$ by $S^1-$valued vector fields
with uniformly bounded $H^1-$norm, it is
necessary for $m'_k$ to have zero degree on each cell of the square
grid ${\cal R}_k$. This property of vanishing degree is shown in
the following lemma:

\begin{lem}
\label{lemcell} Let $0<\alpha<1$, $0<\beta< {1-\alpha}$
and $C>0$. There exists $\eps_0=\eps_0(\alpha, \beta, C)>0$ such that for every $\eps\in (0, \eps_0)$ the following holds:
if ${\cal
Z}=(-\frac{\eps^{\beta}}{2},\frac{\eps^{\beta}}{2})^2$ is the cell
of length $\eps^{\beta}$ and $m':\overline{\cal Z}\to \overline{B^2}$ is a $C^1$ vector field such that 
$$\int_{\partial {\cal Z}} g_\eps(m')\, d\h^1\leq
\frac{C|\log \eps|}{\eps^{\beta}}\quad \textrm{ and } \quad
\int_{{\cal Z}} g_\eps(m')\, dx\leq {2\pi \alpha}|\log \eps|+C\log |\log \eps|,$$
then $\degr(m',
\partial {\cal Z})=0$.
\end{lem}

\proof{ of Lemma \ref{lemcell} } Note that the same argument as in \eqref{modmic} implies that $|m'|\geq 1/2$ on $\partial {\cal
Z}$, so that it makes sense to speak about the degree of $m'$ on $\partial {\cal Z}$. Note also that the quantity $C\log |\log \eps|$ in the upper bound of the GL energy on $\cal Z$ can be absorbed by the leading order term $2\pi \tilde \alpha |\log \eps|$ for a slightly bigger $\tilde \alpha>\alpha$ so that the inequality $\beta<1-\tilde \alpha$ still holds. Therefore, we omit that second leading order term in the following. The idea of the proof consists in a rescaling and extension argument so that the imposed upper bounds on the GL energy rule out the existence of a vortex in the interior. Indeed, assume by contradiction that
$|\degr(m',
\partial {\cal Z})|>0$ (i.e., a vortex exists in the interior). By a change of scale, we define $\tilde m'$
on the rescaled cell ${\cal
Z}_{1/2}=(-1/2, 1/2)^2$ as:
$$\tilde m'(x)=m'(\eps^\beta x)\quad \textrm{ if } x\in {\cal
Z}_{1/2}=(-1/2, 1/2)^2$$ and then, we extend $\tilde
m'$ to the larger cell
${\cal Z}_\lambda=(-\lambda, \lambda)^2$ with
$\lam>1/2$ (to be chosen later) as follows:
$$\tilde m'(x)=\tilde m'(y)\quad \textrm{ if } x\in {\cal Z}_\lambda \setminus
{\cal Z}_{1/2},$$ where $y\in \partial {\cal Z}_{1/2}$ with $y=tx$
for some $t\in (0,1)$ (i.e., $y$ is the closest point to $x$ on
the boundary $\partial {\cal Z}_{1/2}$ that has the same direction
as $x$). Therefore, $|\tilde m'|\geq 1/2$ in ${\cal Z}_\lambda
\setminus {\cal Z}_{1/2}$ and $|\degr(\tilde m', \partial {\cal
Z}_\lambda)|>0$. Letting $\delta=\eps^{1-\beta}$, we have that
$$\int_{{\cal Z}_{1/2}} g_\delta(\tilde m')\, dx=\int_{{\cal Z}} g_\eps(m')\,
dx\leq \frac{2\pi \alpha}{1-\beta}|\log \delta|$$ and
\begin{align*}
\int_{{\cal Z}_\lambda \setminus {\cal Z}_{1/2}} g_\delta(\tilde
m')\, dx&\leq \tilde C (\lambda-\frac 1 2) \int_{\partial
{\cal Z}_{1/2}} g_\delta(\tilde m')\, d\h^1\\
&=\tilde C \eps^\beta (\lambda-\frac 1 2) \int_{\partial {\cal Z}}
g_\eps(m')\, d\h^1\\&\leq \tilde C C (\lambda-\frac 1 2)|\log \eps|=\frac{\tilde C C
(\lambda-\frac 1 2)}{1-\beta}|\log \delta|,
\end{align*}
for some universal constant $\tilde C>0$ and a small $\eps>0$.
We choose $\lambda>1/2$ such that
$$K:=\frac{\alpha}{1-\beta}+\frac{\tilde C C
(\lambda-\frac 1 2)}{2\pi(1-\beta)}<1$$ (this is possible since by hypothesis, $\beta<1-\alpha$). By summing over the above energy estimates,
we obtain that \be \label{est45} \int_{{\cal Z}_{\lambda}}
g_\delta(\tilde m')\, dx\leq 2\pi K|\log
\delta|.\ee Since $K<1$, Theorem \ref{thm_Jer}
implies the existence of a ball $\tilde B\subset {\cal Z}_{\lambda}$ of radius $\lambda-1/2$ with
$$\int_{\tilde B} g_\delta(\tilde m')\, dx\geq 2\pi|\log \delta|-\bar
C$$ for $\delta$ sufficiently small, which is a contradiction with
\eqref{est45}. \qed

\medskip

As a consequence of Lemma \ref{lemcell}, we deduce by \eqref{en_D} and \eqref{condR} that our choice
$\beta< {1-\alpha}$ implies that $m'_k$ has
vanishing degree on every cell of the grid ${\cal R}_k$.

\noindent {\bf Step 4}. {\it Construction of an approximating sequence.} We denote
$$\rho_k=|m'_k| \quad \textrm{ and } \quad m'_k=\rho_k v_k.$$
By {\it Step 3}, we can smoothly lift $m'_k$ on the grid, i.e.,
$$v_k=\frac{m'_k}{\rho_k}=e^{i\vfi_k}\quad \textrm{ on } {\cal R}_k \textrm{   and  } \varphi_k\in C^1({\cal R}_k, \RR).$$
On each cell ${\cal Z}^k$ of length $\eps_k^\beta$ of the grid, we
define
$$M'_k=e^{i\Phi_k}\quad \textrm{ in } {\cal Z}^k$$
where $\Phi_k$ is the harmonic extension of $\vfi_k$ inside ${\cal
Z}^k$, i.e.,
$$\begin{cases}
\Delta \Phi_k=0 & \quad \textrm{ in } {\cal Z}^k\\
\Phi_k=\vfi_k & \quad \textrm{ on } \partial {\cal Z}^k.
\end{cases}
$$
Since $\vfi_k$ can be smoothly extended around $\partial {\cal Z}^k$ (because $m'_k/\rho_k$ has a $C^1$ lifting around $\partial {\cal Z}^k$), we deduce that $\Phi_k\in C^1(\bar{\cal Z}^k)$.
Note that the following inequality holds:
\be
\label{esti_harmon}
\int_{ {\cal
Z}^k} |\nabla
\Phi_k|^2\, dx
\leq C \eps_k^\beta \int_{\partial {\cal Z}^k} |\nabla
\vfi_k|^2\,
d\h^1.\ee
Indeed, after rescaling by $\eps^\beta$, we show the inequality in the unit cell ${\cal Z}_1=(-1,1)^2$ for the harmonic function $\Phi$ in ${\cal Z}_1$ with the trace $\f$ on $\partial {\cal Z}_1$. 
We can assume that $\int_{\partial {\cal Z}_1} \f \,
d\h^1=0$ (otherwise, consider $\f-\int_{\partial {\cal
Z}_1}\hspace{-7.5mm}- \,  \, \,\,\,\,\, \f \,
d\h^1$). For that, we consider
a smooth cut-off function $\Psi:[0,1]\to \RR$ such that $\Psi(t)=0$ for $t\leq 1/2$ and $\Psi(1)=1$ and the following extension $\Phi^{ext}$ of $\f$ in ${\cal Z}_1$: 
$\Phi^{ext}(t\cdot x)=\Psi(t)\f(x)$ for every $t\in (0,1)$ and $x\in \partial {\cal Z}_1$. By Poincar\'e's inequality, one has
$$\int_{ {\cal
Z}_1} |\nabla
\Phi|^2\, dx
\leq \int_{ {\cal
Z}_1} |\nabla
\Phi^{ext}|^2\, dx\leq C \int_{\partial {\cal Z}_1} (|\nabla
\f|^2+ \f^2)\,
d\h^1\leq C \int_{\partial {\cal Z}_1} |\nabla
\f|^2 \,
d\h^1. $$   

The goal is to prove that the sequence $\{M'_k\}$ approximates
$\{m'_k\}$ in $L^2(B^2, \RR^2)$ and $M'_k$ satisfies \eqref{bun_stray} for some associated stray field $h_k$ defined in $B^3$. 

\noindent {\bf Step 5}. {\it Estimate
$\|\nabla(M'_k-m'_k)\|_{L^2}$.}
Denoting by $C$ a generic universal
constant, we have
\begin{align}
\nonumber \int_{ {\cal Z}^k} |\nabla M'_k|^2\, dx&=\int_{ {\cal
Z}^k} |\nabla
\Phi_k|^2\, dx\\
\nonumber&\stackrel{\eqref{esti_harmon}}{\leq} C \eps_k^\beta \int_{\partial {\cal Z}^k} |\nabla
\vfi_k|^2\,
d\h^1\\
\nonumber &= C \eps_k^\beta \int_{\partial {\cal Z}^k} |\nabla
v_k|^2\,
d\h^1\\
\label{cond_grad} &\leq C \eps_k^\beta \int_{\partial {\cal Z}^k}
\rho_k^2 |\nabla v_k|^2\, d\h^1\leq C \eps_k^\beta \int_{\partial
{\cal Z}^k} |\nabla m'_k|^2\, d\h^1\end{align} since $\rho_k\geq
1/2$ on ${\cal R}_k$. Summing up after all cells $ {\cal Z}^k$ of $ {\cal R}_k$, since the convex
hull of ${\cal R}_k$ covers $B^2$, we obtain by \eqref{condR}, \be
\label{desir1} \int_{B^2} |\nabla M'_k|^2\, dx\leq C \eps_k^\beta
\int_{{\cal R}_k} g_{\eps_k}(m'_k) \, d\h^1\leq \frac{C}{\eta_k
|\log \eta_k|}. \ee Combining with \eqref{en_D}, it yields \be
\label{cond-dif_grad} \int_{B^2} |\nabla (M'_k-m'_k)|^2\, dx\leq
\frac{C}{\eta_k |\log \eta_k|}. \ee

\noindent {\bf Step 6}. {\it Estimate $\|M'_k-m'_k\|_{L^2}$.}
By
Poincar\'e's inequality, we have for each cell ${\cal Z}^k$ of
${\cal R}_k$: \be \label{eve1} \int_{ {\cal Z}^k}
\bigg|M'_k-\int_{\partial {\cal Z}^k}\hspace{-7.5mm}- \,\,\,\,\,M'_k
\bigg|^2\, dx\leq C \eps_k^{2\beta} \int_{ {\cal Z}^k} |\nabla
M'_k|^2\, dx\stackrel{\eqref{cond_grad}}{\leq}  C \eps_k^{3\beta}
\int_{\partial {\cal Z}^k} |\nabla m'_k|^2\, d\h^1 \ee and \be
\label{eve2} \int_{ {\cal Z}^k} \bigg|m'_k-\int_{\partial {\cal
Z}^k}\hspace{-7.5mm}- \,\,\,\,\,m'_k \bigg|^2\, dx\leq C
\eps_k^{2\beta} \int_{ {\cal Z}^k} |\nabla m'_k|^2\, dx. \ee Since
$v_k=M'_k$ on $\partial {\cal Z}^k$, by Jensen's inequality, we
also compute
\begin{align}
\nonumber\int_{ {\cal Z}^k} \bigg|\int_{\partial {\cal
Z}^k}\hspace{-7.5mm}- \,\,\,\,\,(M'_k-m'_k) \bigg|^2\, dx&=\int_{
{\cal Z}^k} \bigg|\int_{\partial {\cal Z}^k}\hspace{-7.5mm}-
\,\,\,\,\,(v_k-m'_k) \bigg|^2\, dx\\
\nonumber&\leq C \eps_k^{2\beta} \int_{\partial {\cal
Z}^k}\hspace{-7.5mm}-
\,\,\,\,\,(1-\rho_k)^2\, d\h^1\\
\nonumber&\leq C \eps_k^{\beta} \int_{\partial {\cal Z}^k}(1-\rho_k^2)^2\, d\h^1\\
\label{eve3}&{\leq}  C \eps_k^{\beta+2} \int_{\partial {\cal Z}^k}
g_{\eps_k}(m'_k)\, d\h^1.
\end{align}
Summing up \eqref{eve1}, \eqref{eve2} and \eqref{eve3} over all
the cells ${\cal Z}^k$ of the grid ${\cal R}_k$, by \eqref{en_D} and \eqref{condR}, we
obtain that \be \label{cond-L2} \int_{B^2} |M'_k-m'_k|^2\, dx\leq
\frac{C\eps_k^{2\beta}}{\eta_k |\log \eta_k|}. \ee
\noindent {\bf Step 7}. {\it Construction of an appropriate stray field $h_k$ associated to $M'_k$ in $B^3$ such that \eqref{bun_stray} holds for the couple $(M'_k, h_k)$.} The choice of the stray field $h_k$ has the form
$$h_k:=(\nabla,\frac{\partial}{\partial z}) (U_k+\tilde U_k)$$ where
$U_k$ is the stray field potential associated to $m'_k$ by \eqref{choiceUk} and
we consider $\tilde U_k\in H^1_0(B^3)$ to be the unique solution of the
variational problem
\be \label{relmhu}
\int_{B^3} \left(\nabla \tilde U_k\cdot \nabla
\zeta+\frac{\partial \tilde U_k}{\partial z}\frac{\partial
\zeta}{\partial z}\right)\, dxdz=\int_{B^2} \zeta \nabla\cdot
(M'_k-m'_k)\, dx,\, \forall \zeta\in H^1_0(B^3).
\ee
(It is a direct consequence of Lax-Milgram's Theorem in $H^1_0(B^3)$.) Note that $h_k$ is indeed a stray field associated to $M'_k:B^2\to S^1$ on the unit ball. In order to estimate
$\int_{B^3} |h_k|^2\, dxdz$, we observe that 
\be
\label{starul1}
\int_{B^3} |(\nabla, \frac{\partial}{\partial z}) U_k|^2\, dxdz\leq \int_{\RR^3} |(\nabla, \frac{\partial}{\partial z}) U_k|^2\, dxdz \leq \frac{C}{|\log \eta_k|}\ee by \eqref{choiceUk} and \eqref{en_D}. It remains to estimate $\int_{B^3} |(\nabla, \frac{\partial}{\partial z}) \tilde U_k|^2\, dxdz$. For that, one should use an interpolation argument via \eqref{cond-dif_grad} and \eqref{cond-L2}.
For that, we extend $\tilde U_k$ by $0$ outside $B^3$, so that the extended function (still denoted by $\tilde U_k$) belongs to $H^1(\RR^3)$ and the trace $\tilde U_k\bigg|_{\RR^2}\in H^{1/2}(\RR^2)$. Moreover, we have
\be
\label{pe4}
\int_{\RR^2} ||\nabla|^{1/2} \tilde U_k|^2\, dx\leq \frac 1 2 \int_{\RR^3} \left(|\nabla \tilde U_k|^2+\bigg|\frac{\partial
\tilde U_k}{\partial z}\bigg|^2\right)\, dxdz= \frac 1 2 \int_{B^3}
\left(|\nabla \tilde U_k|^2+\bigg|\frac{\partial
\tilde U_k}{\partial z}\bigg|^2\right)\, dxdz.
\ee
Let us denote by $T$ a linear continuous extension operator:
$$T:H^s(B^2)\to H^s(\RR^2), \, \, s=0,1.$$
Then by interpolation, it follows that
\begin{align*}
\int_{\RR^2} ||\nabla|^{1/2} T(M'_k-m'_k)|^2\, dx&\leq \bigg(\int_{\RR^2} |T(M'_k-m'_k)|^2\, dx\bigg)^{1/2} \bigg(\int_{\RR^2} |\nabla T(M'_k-m'_k)|^2\, dx\bigg)^{1/2}\\
&\leq C\bigg(\int_{B^2} |M'_k-m'_k|^2\, dx\bigg)^{1/2} \bigg(\int_{B^2} |\nabla (M'_k-m'_k)|^2\, dx\bigg)^{1/2}.
\end{align*}
Combining with \eqref{pe4}, the choice $\zeta:=\tilde U_k$ in \eqref{relmhu} yields
\begin{align*}
\nonumber
\int_{B^3}
\left(\big|\nabla \tilde U_k\big|^2+\big|\frac{\partial \tilde
U_k}{\partial z}\big|^2 \right)\, dxdz&=\int_{B^2} \tilde U_k \nabla\cdot
(M'_k-m'_k)\, dx\\
\nonumber
&=\int_{\RR^2} \tilde U_k \nabla\cdot
T(M'_k-m'_k)\, dx\\
&\leq \left(\int_{\RR^2} ||\nabla|^{1/2} \tilde U_k|^2\, dx\right)^{1/2}\left(\int_{\RR^2} ||\nabla|^{1/2} T(M'_k-m'_k)|^2\, dx\right)^{1/2}\\
\nonumber
&\leq C \left(\int_{B^3} |\nabla \tilde U_k|^2+\bigg|\frac{\partial
\tilde U_k}{\partial z}\bigg|^2\, dx dz\right)^{1/2} \\&\quad \times\left(\int_{B^2}
|M'_k-m'_k|^2\, dx\right)^{1/4} \left(\int_{B^2} |\nabla
(M'_k-m'_k)|^2\, dx\right)^{1/4}. \end{align*}
Hence,
\be
\label{star2}
\int_{B^3}
|(\nabla, \frac{\partial}{\partial z}) \tilde U_k|^2\, dxdz
\stackrel{\eqref{cond-dif_grad},
\eqref{cond-L2}}{\leq} \frac{C\eps_k^{\beta}}{\eta_k |\log
\eta_k|}\stackrel{\eqref{new_regi}}{\leq} \frac{C}{|\log
\eta_k|}
\ee
for $k$ sufficiently large. Therefore, by \eqref{starul1} and \eqref{star2}, we conclude
\begin{align}
\label{desir2}
\int_{B^3} |h_k|^2\, dx dz&\leq 
 \frac{C}{|\log \eta_k|}.\end{align} By
\eqref{desir1} and \eqref{desir2}, condition \eqref{bun_stray} is satisfied for $M'_k$ and the stray fields $h_k$. Then Theorem \ref{thm1} applies and implies that $\{M'_k\}$ is
relatively compact in $L^1(B^2)$. Therefore, from \eqref{cond-L2},
it follows that $\{m'_k\}$ also is relatively compact in
$L^1(B^2)$. Since the ball $\bbbb$ was arbitrary chosen in the
complementary of ${\cal D}$ and we proved the relatively
compactness result in the reduced ball $B^2=\frac 1 2 {\bbbb}$, by a diagonal
argument, we deduce that $\{m'_k\}$ converges in
$L^1(\Omega\setminus {\cal D})$ up to a subsequence.
Letting now $\sigma\to 0$, the conclusion of Theorem
\ref{thm_main} follows. \qed

\section{Upper bound for the Landau state}
\label{section_up}

In this section we prove the upper bound stated in Theorem
\ref{thm_upper} for a stadium domain:

\proof{ of Theorem \ref{thm_upper} } The construction is carried out in several steps:

{\noindent {\bf Step 1:} {\it A N\'eel wall approximation}.
Let $$\lambda:={\eta}{|\log
\eta|}.$$ The parameter $\lambda$ corresponds to the core size of a $180^\circ$ wall transition. More precisely, we consider the following $1d$ transition
layer $(u_\lambda, v_\lambda):\RR\to S^1$ that approximates a $180^\circ$ N\'eel
wall centered at the origin (see Figure \ref{neelu}): 
$$u_\lambda(t)=\begin{cases} \frac{|\log \sqrt{t^2+\lambda^2}|}{|\log
\lambda|}&\quad \textrm{ if } |t|\leq \sqrt{1-\lambda^2},\\
0 &\quad \textrm{ elsewhere},
\end{cases}$$
$$v_\lambda(t)=\begin{cases} -\sqrt{1-u^2_\lambda(t)}
&\quad \textrm{ if } t\leq 0,\\
\sqrt{1-u^2_\lambda(t)} &\quad \textrm{ if } t\geq 0.
\end{cases}
$$

\begin{figure}[htbp]
\center
\includegraphics[scale=0.3,
width=0.3\textwidth]{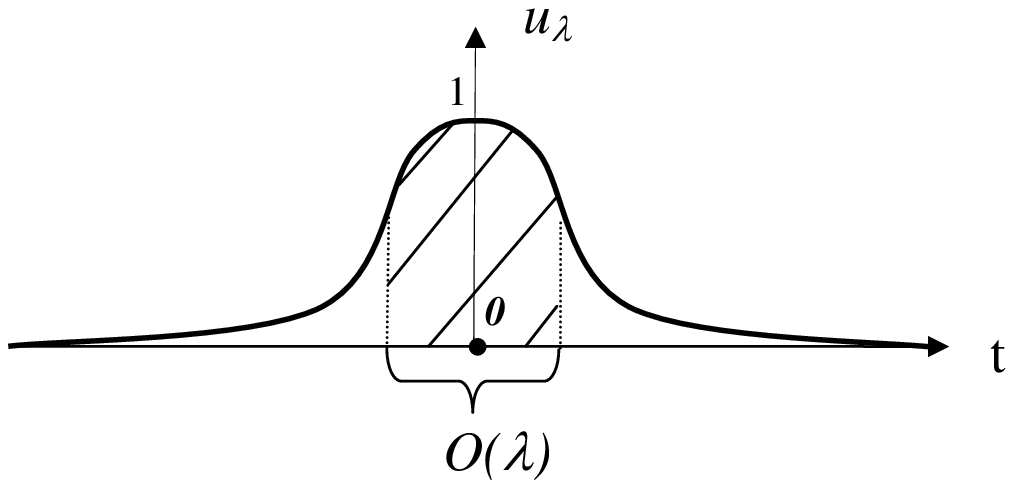} 
\includegraphics[scale=0.3,
width=0.3\textwidth]{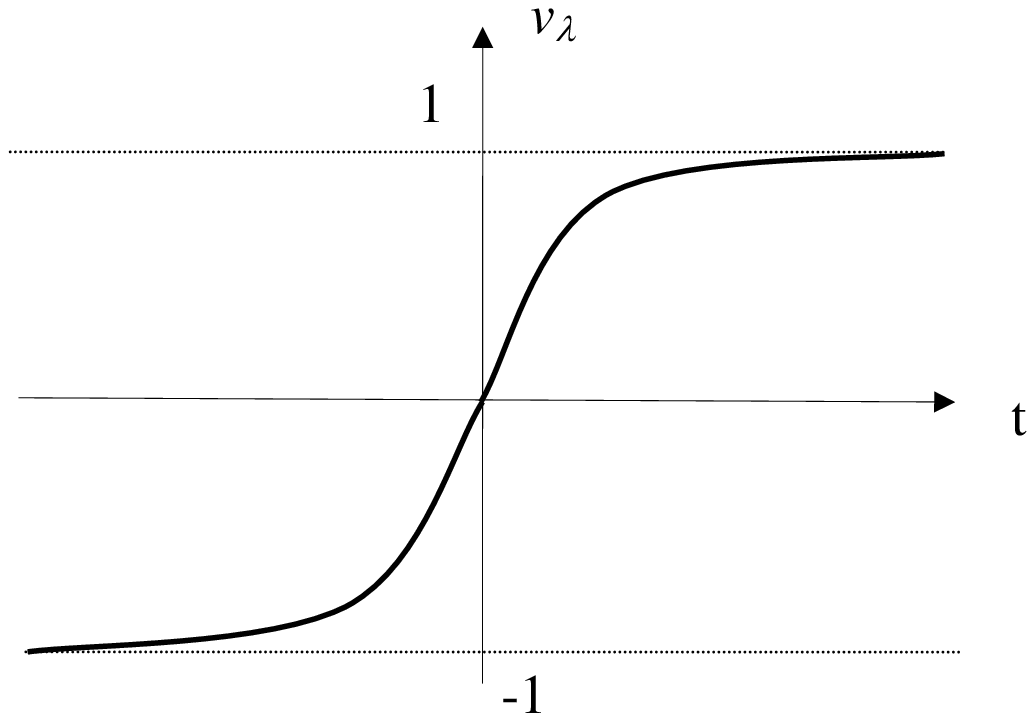} 
\caption{A $180^\circ$ N\'eel
wall approximation.}
\label{neelu}
\end{figure}

The exchange energy corresponding to this transition layer
estimates as follows (see DeSimone, Kn\"upfer and Otto
\cite{DKO} or Ignat \cite{IGNAT_gamma}): \be \label{Extran}\int_\RR
\big|\frac{du_\lambda}{dt}\big|^2+\big|\frac{dv_\lambda}{dt}\big|^2\, dt
\leq \int_\RR \frac{1}{1-u_\lambda} \big|\frac{du_\lambda}{dt}\big|^2 \,
dt=O\bigg(\frac{1}{\lambda |\log \lambda|}\bigg). \ee In order to
estimate the stray-field energy of the transition layer, let
$U_\lambda$ be the radial extension of $u_\lambda$ in $\RR^2$:
$$U_\lambda(x_1, x_2)=u_\lambda(\sqrt{x_1^2+x_2^2}\, ).$$
By $\hd-$trace estimate of an $\hu-$function, it follows (see details in
\cite{DKO, IGNAT_gamma} or \eqref{esti_h12_ex} below): 
\be \label{Straytran}
\|u_\lambda\|_{\hd}^2\leq \frac 1 2 \int_{\RR^2} \big|\nabla
U_\lambda|^2\, dx\leq \pi \int_{0}^1 r
\big|\frac{du_\lambda}{dr}\big|^2\,
dr=\frac{\pi+o(1)}{|\log \eta|}.
\ee

We will construct a continuous vector field $m:\Omega\to S^2$ such that the upper bound in Theorem~\ref{thm_upper} holds and
\be
\label{askin}
m'(x)=\nu^\perp(x), \,\, m_3(x)=0 \quad \textrm{ if }\quad x\in \partial \Omega,\ee
where $\nu$ is the outer unit normal vector on $\partial \Omega$. Moreover, the function $m$ will satisfy the following symmetry properties:
$$m'(x)=-m'(-x), \, \, m_1(x)=-m_1(x_1, -x_2),  \, \, m_2(x)=-m_2(-x_1, x_2), \quad x\in \Omega.$$

 \bigskip

{\noindent {\bf Step 2:} {\it Construction in $\Omega_1$ (the sub-domain defined in Theorem \ref{thm_upper}) }.} We distinguish two regions in $\Omega_1$ (see Figure \ref{ome1}): $$\Omega_{1,1}=\{x\in
\Omega_1\, : \, x_1\geq 1+\delta\}\quad \textrm{ and }\quad
\Omega_{1,2}=\{x\in \Omega_1\, : \, 1\leq x_1< 1+\delta\} \quad \textrm{ with }\quad \delta=\frac{1}{|\log \eta|^{3/2}}.$$ 

\begin{figure}[htbp]
\center
\includegraphics[scale=0.1,
width=0.1\textwidth]{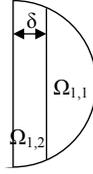} 
\caption{The region $\Omega_1$.} \label{ome1}
\end{figure}

In
$\Omega_{1,1}$, we define $m$ with values in $S^1$ that behaves like a
vortex centered in $A=(1,0)$:
$$m'(x)=\left(\frac{x-A}{|x-A|}\right)^\perp,\, \, m_3(x)=0 \quad \textrm{ in } \quad \Omega_{1,1}.$$
By setting $m'$ to be a $180^\circ$ transition wall on $\partial \Omega_2\cap \partial \Omega_{1,2}$ (as in Step 1), i.e.,
$$m'(1, x_2)=(u_\lambda(x_2), v_\lambda(x_2))^\perp=e^{i\t(x_2)}, \, \, m_3=0 \quad \textrm{ if }\quad x_2\in (-1,1),$$
the vector field $m$ is completely defined on $\partial \Omega_{1,2}$ (together with \eqref{askin}).
Here, $\t$ is the angle transition between $[0, \pi]$ of the $180^\circ$ wall on $\partial \Omega_{1,2}\cap \partial \Omega_2$, i.e., 
$\t(x_2)=0$ and $\t(-x_2)=\pi$ if $x_2\in [-1, -\sqrt{1-\lambda^2}]$, and
$$\t(x_2)=\arcsin\bigg( \frac{1}{|\log \lambda|} \log \frac{1}{\sqrt{x_2^2+\lambda^2}}\bigg),  \, \, \t(-x_2)=\pi-\t(x_2) \quad \textrm{ if }\quad x_2\in [-\sqrt{1-\lambda^2},0].$$
Therefore, we define $m'=e^{i\f}, m_3=0$ inside $\Omega_{1,2}$ by a phase $\f$ that is uniquely determined by the boundary conditions
on $\partial \Omega_{1,2}$ as an affine continuous function in $x_1$:
\begin{align*}
\f(1+t\sqrt{1-x_2^2}, x_2)&=t\arcsin \sqrt{1-x_2^2}+(1-t)\t(x_2), \quad t\in (0,1), \, \,  x_2\in (-1, -\sqrt{1-\delta^2}),\\
\f(1+\delta t, x_2)&=t \arcsin \frac{\delta}{\sqrt{x_2^2+\delta^2}}+(1-t)\t(x_2),\quad  t\in (0,1), \, \, x_2\in (-\sqrt{1-\delta^2},0),\\
\f(x_1, x_2)&=\pi-\f(x_1, -x_2), \quad  x\in \Omega_{1,2},\, \,  x_2>0.\end{align*}
We will denote by $$\alpha_\delta(x_2)=\arcsin \frac{\delta}{\sqrt{x_2^2+\delta^2}}$$ the phase of the vortex at $\partial \Omega_{1,1}\cap \partial \Omega_{1,2}$. 

{\noindent {\bf Step 3:} {\it Estimate of the exchange energy in $\Omega_1$}.} 
First, we have:
\be \label{enex4}
\int_{\Omega_{1,1}} |\nabla m|^2\, dx=O(|\log{\delta}|)=O(\log |\log{\eta}|),\ee
$$
\mathop{\mathop{\int \int}_{x\in \Omega_{1,2}}}_{|x_2|\in (\sqrt{1-\delta^2}, 1)} |\nabla m|^2\, dx=2
\mathop{\mathop{\int \int}_{x\in \Omega_{1,2}}}_{x_2\in (-1, -\sqrt{1-\delta^2})} |\nabla \f|^2\, dx
=o(\delta)
$$
and
\begin{align*}
\mathop{\mathop{\int \int}_{x\in \Omega_{1,2}}}_{|x_2|< \sqrt{1-\delta^2}} |\nabla m|^2\, dx&=2
\mathop{\mathop{\int \int}_{x\in \Omega_{1,2}}}_{x_2\in (-\sqrt{1-\delta^2}, 0)} |\nabla \f|^2\, dx\\
&=2\int_0^1 \int_{-\sqrt{1-\delta^2}}^0 \bigg(\frac{1}{\delta}(\alpha_\delta(x_2)-\t(x_2))^2+\delta(t \frac{d\alpha_\delta}{dx_2}+(1-t)\frac{d\t}{dx_2})^2
\bigg)\, dt dx_2\\
&\leq 4\int_0^{\sqrt{1-\delta^2}} \bigg(\frac{1}{\delta}\alpha_\delta^2(x_2)+\frac 1 \delta \t^2(x_2)+\delta \big|\frac{d\alpha_\delta}{dx_2}\big|^2+
\delta \big|\frac{d\t}{dx_2}\big|^2
\bigg)\, dx_2.
\end{align*}
Introducing the notation $\alpha(\frac{x_2}{\delta})=\alpha_\delta(x_2)$,
we compute:
\be
\label{alf}
\int_0^{\sqrt{1-\delta^2}} \frac{1}{\delta}\alpha_\delta^2(x_2)\, dx_2\leq \int_0^{1/\delta} \alpha^2(s)\, ds\leq 4 \int_0^{1/\delta} \frac{1}{s^2+1}\, ds=O(1)\ee
(where we use that $\arcsin x\leq 2x$ if $x\in (0,1)$),
\be
\label{the}
\int_0^{\sqrt{1-\delta^2}} \frac{1}{\delta}\t^2(x_2)\, dx_2=O(\frac{1}{\delta|\log \lambda|^2})=o(1),\ee
$$\int_0^{\sqrt{1-\delta^2}} {\delta} \big|\frac{d\alpha_\delta}{dx_2}\big|^2 \, dx_2\leq \int_0^{1/\delta} \big|\frac{d\alpha}{ds}\big|^2(s)\, ds= \int_0^{1/\delta} \frac{1}{(s^2+1)^2}\, ds=O(1)$$
and
$$\int_0^{\sqrt{1-\delta^2}} {\delta} \big|\frac{d\t}{dx_2}\big|^2 \, dx_2=O(\frac{\delta}{\lambda|\log \lambda|})=o(\frac{1}{\eta|\log \eta|}).$$
Therefore, \be \label{enex5}
\int_{\Omega_{1,2}} |\nabla m|^2\, dx=o(\frac{1}{\eta|\log \eta|}).\ee

\bigskip

{\noindent {\bf Step 4:} {\it Construction in $\Omega_3$}.} We define $m$ by imposing the symmetry
$m(x)=-m(-x)$ for $x\in \Omega_3$. Therefore, by \eqref{enex4} and \eqref{enex5}, we have
\be \label{enexom3tot}
\int_{\Omega_3} |\nabla m|^2\, dx=\int_{\Omega_1} |\nabla m|^2\, dx=o(\frac{1}{\eta|\log \eta|}).\ee

{\noindent {\bf Step 5:} {\it Construction in $\Omega_2$}.} We distinguish two regions in $\Omega_2$ (see Figure \ref{ome2}): $$\Omega_{2,1}=\{x\in
\Omega_2\, : \, |x_1|\in (2\delta, 1)\}\quad \textrm{ and }\quad
\Omega_{2,2}=\{x\in \Omega_1\, : \, |x_1|< 2\delta\}.$$ 

\begin{figure}[htbp]
\center
\includegraphics[scale=0.3,
width=0.3\textwidth]{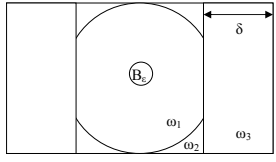} 
\caption{The region $\Omega_{2,2}$.} \label{ome2}
\end{figure}

In
$\Omega_{2,1}$, we define $m$ with values in $S^1$ that behaves like a
$180^\circ$ N\'eel wall (as in Step 1):
\begin{align*}
m'(x)&=(u_\lambda(x_2), v_\lambda(x_2))^\perp=e^{i\t(x_2)},\, \, m_3(x)=0 \quad \textrm{ for } \quad  x_1\in (2\delta, 1), \, \,  x_2\in (-1, 1), \\
m(x)&=-m(-x) \quad \textrm{ for } \quad x_1\in (-1, -2\delta),  \, \,  x_2\in (-1, 1).
\end{align*}
Denoting by $B_r$ the disc centered at the origin of radius $r$, we decompose the domain $$\Omega_{2,2}=B_\eps \cup \omega_1\cup \omega_2\cup \omega_3$$
with
\begin{align*}
\omega_1&=\{x\in B_1\, :\, |x_1|\leq \delta\}\setminus B_\eps,\\
\omega_2&=\{x\in \Omega_2\setminus B_1 \, :\, |x_1|\leq \delta\},\\
\omega_3&=\{\delta<|x_1|<2\delta\}\times (-1,1).
\end{align*}
In $B_\eps$ (the core of the vortex), we define
$$m'(x)=\sin(\frac{\pi}{2\eps}|x|) \left(\frac{x}{|x|}\right)^\perp,\,\, m_3(x)=\sqrt{1-|m'|^2(x)} \quad \textrm{ for } \quad x\in B_\eps. $$
In $\omega_1$, we define $m$ with values in $S^1$ that corresponds to the vortex away from the core:
$$m'(x)=\left(\frac{x}{|x|}\right)^\perp,\,\, m_3(x)=0 \quad \textrm{ for } \quad x\in B_1\setminus B_\eps \textrm{ and } |x_1|\leq \delta.$$
In $\omega_2$, we define $m$ with values in $S^1$: $m'=e^{i\f}, m_3=0$ inside $\omega_2$. The phase $\f$ is given as an affine continuous function in $x_2$ determined by the values on the boundary $\partial \omega_2$:
\begin{align*}
\f(x_1, -(1-t)-t\sqrt{1-x_1^2})&=t\arcsin {x_1}, \quad t\in (0,1), x_1\in (0, \delta),\\
\f(x_1, x_2)&=\pi-\f(x_1, -x_2), \quad  x\in \omega_2, x_1\in (0, \delta), x_2>0,\\
\f(x)&=\pi+\f(-x), \quad  x\in \omega_2, x_1\in(-\delta, 0).
\end{align*}
In $\omega_3$, we also define $m$ with values in $S^1$ where the phase $\f$ is an affine continuous function in $x_1$ determined by the boundary conditions
on $\partial \omega_2$:
\begin{align*}
\f(\delta+\delta t, x_2)&=\frac{(1-t)(x_2+1)}{1-\sqrt{1-\delta^2}} \arcsin {\delta}+t\t(x_2),\quad  t\in (0,1), x_2\in (-1,-\sqrt{1-\delta^2}),\\
\f(\delta+\delta t, x_2)&=(1-t) \alpha_{\delta}(x_2)+t\t(x_2),\quad  t\in (0,1), x_2\in (-\sqrt{1-\delta^2},0),\\
\f(x_1, x_2)&=\pi-\f(x_1, -x_2), \quad  x_1\in (\delta, 2\delta), x_2\in (0, 1),\\
\f(x)&=\pi+\f(-x), \quad  x_1\in (-2\delta, -\delta), x_2\in (-1, 1).
\end{align*}

{\noindent {\bf Step 6:} {\it Estimate of the exchange energy in $\Omega_2$}.} 
We start by estimating the exchange energy in $\Omega_{2,1}$ and then, in $\Omega_{2,2}$. By \eqref{Extran}, we have that
\be \label{enex6}
\int_{\Omega_{2,1}} |\nabla m|^2\, dx=2(1-2\delta)\int_\RR
\bigg(\big|\frac{du_\lambda}{dt}\big|^2+\big|\frac{dv_\lambda}{dt}\big|^2 \, \bigg)\, dt=o\bigg(\frac{1}{\eta|\log \eta|}\bigg). \ee
In $\Omega_{2,2}$, we first have
\be \label{enexDeps}
\int_{B_\eps} |\nabla m|^2\, dx=O(1).\ee
Then
\be \label{enexD1}
\int_{\omega_1} |\nabla m|^2\, dx=2\pi|\log \eps|-O(|\log \delta|),\ee
\be \label{enexom2}
\int_{\omega_2} |\nabla m|^2\, dx=\int_{\omega_2} |\nabla \f|^2\, dx=o(\delta)\ee
and
\be \label{enexom3}
\int_{\omega_3} |\nabla m|^2\, dx=\int_{\omega_3} |\nabla \f|^2\, dx\stackrel{\eqref{enex5}}{=}o(\frac{1}{\eta|\log \eta|}).\ee
By \eqref{enexom3tot},  \eqref{enex6},  \eqref{enexDeps},  \eqref{enexD1},  \eqref{enexom2} and  \eqref{enexom3}, we deduce the following estimate of the exchange energy of $m$:
 \be \label{enextotal}
\int_{\Omega} |\nabla m|^2\, dx=
2\pi |\log \eps|+o\bigg(\frac{1}{\eta|\log \eta|}\bigg). \ee

\bigskip

{\noindent {\bf Step 7:} {\it Symmetries of the stray field}.}
Now we estimate the stray field energy of $m$. For that, let $U\in {\cal BL}$ be the stray field potential in the Beppo-Levi space associated to
$m$ defined by \eqref{1c} for $(\nabla,\frac{\partial}{\partial
z}) U$ that satisfies
\be
\label{choiceUuu}
\int_{\RR^2\times \RR}\left(|\nabla U|^2+\big|\frac{\partial U}{\partial
z}\big|^2 \right)\,dx dz=\int_{\Omega} U(x,0)\nabla\cdot m'(x)\, dx= \frac 1 2 \int_{\RR^2}
\left|\,|\nabla|^{-1/2}(\nabla\cdot m')\right|^2\, dx.\ee
Moreover, the stray field potential verifies:
$$\begin{cases}
\Delta U=0\quad &\textrm{ if } z\neq 0,\\
\bigg[\frac{\partial U}{\partial z}\bigg]=-\nabla\cdot m', \,\,  [U]=0 \quad &\textrm{ if } z=0.
\end{cases} $$ 
Since $m'$ is a Lipschitz vector field in $\Omega$ (so, $\nabla\cdot m'\in L^\infty(\RR^2)$), by standard regularity theory for elliptic PDEs, we know that $U$ is continuous in $\RR^3$.
We also may deduce some symmetry properties of $U$: First of all, the uniqueness of the stray field potential $U\in {\cal BL}$ in \eqref{choiceUuu} yields
$U(x, z)=U(x, -z)$ for every $(x,z)\in \RR^3$.
Also, remark that our vector field $m'$ is anti-symmetric with respect to the origin, i.e., $m'(x)=-m'(-x)$ which yields
$\nabla\cdot m'(x)=\nabla\cdot m'(-x)$ in $\RR^2$. Again, by the uniqueness of the stray field potential $U\in {\cal BL}$,
we deduce that $U$ is symmetric in the in-plane variables with respect to the origin, i.e., $$U(x, z)=U(-x, \pm z) \quad \textrm{ for every }\quad (x,z)\in \RR^3.$$ Also, the vector field $m$ satisfies the symmetry relation $m'(x)=(m_1, -m_2)(-x_1, x_2)$ in $\RR^2$, so that $\nabla\cdot m'(x)=-\nabla\cdot m'(-x_1, x_2)$ in $\RR^2$. It implies that $$U(x_1,x_2, z)=-U(-x_1, x_2, \pm z) \quad \textrm{ for every }\quad (x,z)\in \RR^3.$$
Similarly, $U(x_1,x_2, z)=-U(x_1, -x_2, \pm z)$ for every $(x,z)\in \RR^3$. In particular, it yields $U(0, x_2, z)=U(x_1,0,z)=0$ for every $(x,z)\in \RR^3$.

In what follows, we compute upper bounds for  $\ds \int_{\Omega} U(x,0)\nabla\cdot m'(x)\, dx$ in several steps corresponding to each subdomain
of $\Omega$. In $\Omega_{1,1}\cup (-\Omega_{1,1})\cup \omega_1$, $m'$ is of vanishing divergence, therefore \be \label{stra4}
\int_{\Omega_{1,1}\cup (-\Omega_{1,1})\cup \omega_1} U(x,0)\nabla\cdot m'(x)\, dx=0.\ee
In the next step, we estimate
\be
\label{neast}
\int_{\tilde \omega} U(x,0)\nabla\cdot m'(x)\, dx\ee
where $$\tilde \omega=\Omega_{1,2} \cup (-\Omega_{1,2}) \cup \Omega_{2,1}\cup \omega_2\cup \omega_3.$$
In the last step, we compute $\int_{B_\eps} U(x,0)\nabla\cdot m'(x)\, dx$.

\bigskip

{\noindent {\bf Step 8:} {\it Upper bound for \eqref{neast}.} The computation will be done according to the decomposition: $\ds \nabla\cdot m'=\frac{\partial m_1}{\partial x_1}+\frac{\partial m_2}{\partial x_2}$.
In order to estimate $\ds \int_{\Omega_{1,2}\cup (-\Omega_{1,2})} U(x,0)\frac{\partial m_2}{\partial x_2}(x)\, dx$, we use the following argument (see also Proposition 3 in \cite{IGNAT_gamma}):
\begin{lem}
\label{lemh12}
Let $L>0$, $U:\RR^2\to \RR$ and $v:\RR\to \RR$ be such that $v(x_1)=v(-L)=v(L)$ for every $|x_1|\geq L$.
Then $$\bigg(\int_{-L}^L U(x_1,0)\frac{\partial v}{\partial x_1}(x_1)\, dx_1\bigg)^2\leq \frac 1 2 \|v\|^2_{\hd}\bigg(\int_{\RR^2}|\nabla U|^2\, dx\bigg),$$
where
\be
\label{esti_h12_ex} \|v\|^2_{\hd}=\frac 1 2 \min \{\int_{\RR^2}|\nabla V|^2\, dx\, :\, V(x_1, 0)=v(x_1) \textrm{ for every }x_1\in \RR\}.  \ee
\end{lem}

Here, we denote the homogeneous $\dot{H}^{1/2}$-seminorm of $v$ by
$$\|v\|_{\hd}:=\int_\RR |\xi| |\caf{v}|^2(\xi)\, d\xi,$$
where  $\caf{v}\in {\cal S}'(\RR)$ stands for the Fourier
transform of $v$ (as a tempered distribution), i.e.,
$$\caf{v}(\xi)=\frac{1}{\sqrt{2\pi}}\int_\RR e^{-i\xi x}v(x_1)\, dx_1, \quad \forall \xi \in \RR.$$
One can also write
\be
\label{caracte}
\|v\|_{\hd}^2=\frac{1}{2\pi} \int_\RR \int_\RR \frac{|v(s)-v(t)|^2}{|s-t|^2}\, ds dt\ee
(see e.g., \cite{IGNAT_gamma}). Another remark is that for even functions $v$ (i.e., $v(x_1)=v(-x_1)$), the following estimate of $\|v\|_{\hd}$ can be obtained via
\eqref{esti_h12_ex} by
considering the radial extension $V$ of $v$ in $\RR^2$ (i.e., $V(x)=v(|x|)$):
\be
\label{esti_h12_ra}
\|v\|^2_{\hd}\leq \frac 1 2 \int_{\RR^2}|\nabla V|^2\, dx=\pi \int_{0}^L r\big| \frac{\partial v}{\partial r}\big|^2\, dr.
\ee
Observe that \eqref{esti_h12_ex} is a general characterization of the $H^{1/2}-$trace of $H^1-$functions and it is valid in any dimension.

\proof{ of Lemma \ref{lemh12} } W.l.o.g., we can assume that $v(x_1)=v(-L)=v(L)=0$ for every $|x_1|>L$. Then Parseval's identity and the Cauchy-Schwarz inequality yield:
\begin{align*}
\bigg(\int_{-L}^L U(x_1,0)\frac{\partial v}{\partial x_1}(x_1)\, dx_1\bigg)^2&=\bigg(\int_{\RR} U(x_1,0)\frac{\partial v}{\partial x_1}(x_1)\, dx_1\bigg)^2\\
&=\bigg(\int_{\RR} \caf(U(\cdot,0))(\xi_1) \overline{\caf(\frac{\partial v}{\partial x_1})(\xi_1)}\, d\xi_1\bigg)^2\\
&=\bigg(\int_{\RR} i\xi_1 \caf(U(\cdot,0))(\xi_1) \overline{\caf(v)(\xi_1)}\, d\xi_1\bigg)^2\\
&\leq \bigg(\int_{\RR} |\xi_1| \, |\caf(U(\cdot,0))(\xi_1)|^2 \, d\xi_1 \bigg) \bigg(\int_{\RR} |\xi_1| \, |\caf(v)(\xi_1)|^2 \, d\xi_1\bigg)\\
&\leq \|v\|^2_{\hd} \|U(\cdot,0)\|^2_{\hd}\\
&\stackrel{\eqref{esti_h12_ex}}{\leq} \frac 1 2 \|v\|^2_{\hd}\bigg(\int_{\RR^2}|\nabla U|^2\, dx\bigg).
\end{align*} \qed

\bigskip

Writing each $x_1\in (1, 1+\delta)$ as $x_1=1+\delta t$ with $t\in (0,1)$, the $x_2-$section in $\Omega_{1,2}$ passing through $x_1$ is given by
$$\Ix=(-\sqrt{1-\delta^2 t^2}, \sqrt{1-\delta^2 t^2}),$$
so that $\Omega_{1,2}=\cup_{t\in (0,1)} \{x_1\}\times \Ix$. Since $m_2(x_1, \cdot)$ takes the same value at the boundary $\partial \Ix$ for every $t\in (0,1)$, we have by
\eqref{esti_h12_ra} that:
$$\|m_2(x_1, \cdot)\|^2_{\hd}\leq \frac \pi 2 \int_{\Ix} |x_2|\big| \frac{\partial m_2}{\partial x_2}(x_1, x_2)\big|^2\, dx_2\leq \frac \pi 2 \int_{\Ix} |x_2|\big| \frac{\partial \f}{\partial x_2}(x_1, x_2)\big|^2\, dx_2=O(1),
$$
where the upper bound $O(1)$ does not depend on $x_1$. Therefore,
Lemma \ref{lemh12} yields:
\begin{align}
\nonumber
\int_{\Omega_{1,2} \cup (-\Omega_{1,2})} U(x,0)\frac{\partial m_2}{\partial x_2}(x)\, dx &\stackrel{x_1=1+t\delta}{=}2\int_0^1\delta \bigg(\int_{\Ix} U(x, 0) \frac{\partial m_2}{\partial x_2}(x)\, dx_2\bigg)\, dt\\
\nonumber&\leq \sqrt{2} \int_0^1\delta \bigg(\int_{\RR^2} |\bigg(\frac{\partial}{\partial x_2}, \frac{\partial}{\partial z} \bigg)U(1+t\delta, x_2, z)|^2 \, dx_2\, dz\bigg)^{1/2} \|m_2(x_1, \cdot)\|_{\hd}  \, dt\\
\label{stra5_2}
&\leq C \sqrt{\delta} \bigg(\int_{\RR^3} |\bigg(\nabla, \frac{\partial}{\partial z} \bigg) U(x, z)|^2 \, dx\, dz\bigg)^{1/2}.
\end{align}
We apply the same argument to estimate $\ds \int_{\Omega_{2,1}} U(x,0)\frac{\partial m_2}{\partial x_2}(x)\, dx$. By \eqref{Straytran}, we already know that
$$\|m_2(x_1, \cdot)\|^2_{\hd}=\frac{\pi+o(1)}{|\log
\eta|}, \textrm{ for all } \, |x_1|\in (2\delta, 1).$$
We deduce via Lemma \ref{lemh12} that:
\begin{align}
\nonumber
\int_{\Omega_{2,1}} U(x,0)\frac{\partial m_2}{\partial x_2}(x)\, dx &{=}\int_{2\delta <|x_1|<1} \bigg(\int_{-1}^1 U(x, 0) \frac{\partial m_2}{\partial x_2}(x)\, dx_2\bigg)\, dx_1\\
\nonumber&\leq \frac{1}{\sqrt{2}} \int_{2\delta <|x_1|<1}  \bigg(\int_{\RR^2} |\bigg(\frac{\partial}{\partial x_2}, \frac{\partial}{\partial z} \bigg)U(x, z)|^2 \, dx_2\, dz\bigg)^{1/2} \|m_2(x_1, \cdot)\|_{\hd}  \, dx_1\\
\label{stra6_2}
&\leq \left( \frac{\pi+o(1)}{|\log
\eta|}\right)^{1/2}\bigg(\int_{\RR^3} |\bigg(\nabla, \frac{\partial}{\partial z} \bigg)U(x, z)|^2 \, dx\, dz\bigg)^{1/2}.
\end{align}
When estimating the same quantity in $\omega_3$, a similar computation to \eqref{stra5_2} leads to
\be
\label{strayom3_2}
\int_{\omega_3} U(x,0)\frac{\partial m_2}{\partial x_2}(x)\, dx \leq C \sqrt{\delta} \bigg(\int_{\RR^3} |\nabla_{(x,z)} U(x, z)|^2 \, dx\, dz\bigg)^{1/2}.
\ee
In $\omega_2$, a slightly different argument is used to estimate the quantity:
\begin{align}
\nonumber
\int_{\omega_2} U(x,0)\frac{\partial m_2}{\partial x_2}(x)\, dx &\leq \bigg(\int_{\omega_2} |U(x, 0)|^4\, dx\bigg)^{1/4}  \bigg(\int_{\omega_2} \bigg|\frac{\partial m_2}{\partial x_2}\bigg|^{4/3}\, dx\bigg)^{3/4}\\
\nonumber &\leq  \bigg(\int_{\RR^2} |U(x, 0)|^4\, dx\bigg)^{1/4}  \bigg(\int_{\omega_2} \bigg|\frac{\partial \f}{\partial x_2}\bigg|^{4/3}\, dx\bigg)^{3/4}\\
\nonumber &\leq C {\delta}^{5/4} \|U(\cdot, 0)\|_{\dot{H}^{1/2}(\RR^2)}\\
\label{straom2_2}
&\stackrel{\eqref{esti_h12_ex}}{\leq} C {\delta}^{5/4} \bigg(\int_{\RR^3} |\bigg(\nabla, \frac{\partial}{\partial z} \bigg) U(x, z)|^2 \, dx\, dz\bigg)^{1/2}.
\end{align}

It remains to estimate $\ds \int_{\tilde \omega} U(x,0)\frac{\partial m_1}{\partial x_1}(x)\, dx$.
In the region near the boundary, i.e., $\tilde \omega\cap \{\sqrt{1-\delta^2}\leq |x_2|\leq 1\}$, the same argument as in \eqref{straom2_2} yields:
\begin{align}
\nonumber
\int_{\tilde \omega\cap\{\sqrt{1-\delta^2}\leq |x_2|\leq 1\}} U(x,0)\frac{\partial m_1}{\partial x_1}(x)\, dx &\leq  \bigg(\int_{\RR^2} |U(x, 0)|^4\, dx\bigg)^{1/4}
\bigg(\int_{\tilde \omega\cap\{\sqrt{1-\delta^2}\leq |x_2|\leq 1\}} \bigg|\frac{\partial \f}{\partial x_1}\bigg|^{4/3}\, dx\bigg)^{3/4}\\
\label{straomtil_1}
&\leq C {\delta}^{9/4} \bigg(\int_{\RR^3} |\bigg(\nabla, \frac{\partial}{\partial z} \bigg) U(x, z)|^2 \, dx\, dz\bigg)^{1/2}.
\end{align}
For the interior region, i.e., $\tilde \omega\cap \{|x_2|\leq \sqrt{1-\delta^2}\}$, we notice that $\frac{\partial m_1}{\partial x_1}\equiv 0$ on $\Omega_{2,1}$ and 
up to a translation, $\frac{\partial m_1}{\partial x_1}$ coincides on $\Omega_{1,2}$ and $\omega_3$. Therefore, it is enough to estimate (by the above argument) the quantity
\begin{align*}
\int_{\omega_3\cap \{|x_2|\leq \sqrt{1-\delta^2}\}} U(x,0)\frac{\partial m_1}{\partial x_1}(x)\, dx &\leq \bigg(\int_{\RR^2} |U(x, 0)|^4\, dx\bigg)^{1/4}  \bigg(\int_{\omega_3\cap \{|x_2|\leq \sqrt{1-\delta^2}\}} 
\bigg|\frac{\partial \f}{\partial x_1}\bigg|^{4/3}\, dx\bigg)^{3/4}\\
 &\leq \bigg(\int_{\RR^2} |U(x, 0)|^4\, dx\bigg)^{1/4}  \bigg(\int_{0}^{\sqrt{1-\delta^2}} \frac{1}{\delta^{1/3}}  (\alpha_\delta(x_2)-\t(x_2))^{4/3}\, dx_2\bigg)^{3/4}.
\end{align*}
The same computation as in \eqref{alf} and \eqref{the} yields
$$ \frac{1}{\delta^{1/3}} \int_{0}^{\sqrt{1-\delta^2}}  \alpha_\delta^{4/3}(x_2) \, dx_2\lesssim  \delta^{2/3}\int_0^{1/\delta} \frac{1}{(t^2+1)^{2/3}}\, dt=O(\delta^{2/3})$$
and
$$ \frac{1}{\delta^{1/3}} \int_{0}^{\sqrt{1-\delta^2}}  \t^{4/3}(x_2) \, dx_2=O(\frac{1}{\delta^{1/3} |\log \lambda|^{4/3}}).$$
Therefore, we deduce that
\be
\label{stra_nustiu}
\int_{\tilde \omega\cap \{|x_2|\leq \sqrt{1-\delta^2}\}} U(x,0)\frac{\partial m_1}{\partial x_1}(x)\, dx
{\leq} o(\frac{1}{|\log \eta|^{1/2}}) \bigg(\int_{\RR^3} |\bigg(\nabla, \frac{\partial}{\partial z} \bigg) U|^2 \bigg)^{1/2}.
\ee
Summing \eqref{stra5_2}, \eqref{stra6_2}, \eqref{strayom3_2}, \eqref{straom2_2}, \eqref{straomtil_1} and \eqref{stra_nustiu}, we obtain the following estimate for the stray field energy in $\tilde \Omega$:
\be
\label{stratil}
\int_{\tilde \omega} U(x,0) \nabla\cdot m'\, dx {\leq} \bigg(\frac{\pi+o(1)}{|\log
\eta|}\bigg)^{1/2} \bigg(\int_{\RR^3} |\bigg(\nabla, \frac{\partial}{\partial z} \bigg) U|^2 \bigg)^{1/2}.
\ee

\bigskip

{\noindent {\bf Step 9:} {\it Conclusion.}
It remains to estimate the stray field energy in $B_\eps$ as in \eqref{straom2_2}:
\begin{align}
\nonumber
\int_{B_\eps} U(x,0)\nabla \cdot m'\, dx &{\leq} C\bigg(\int_{\RR^3} |\bigg(\nabla, \frac{\partial}{\partial z} \bigg) U|^2 \bigg)^{1/2} \bigg(\int_{B_\eps} |\nabla m'|^{4/3}\bigg)^{3/4}\\
\label{stradep} &{\leq} C\sqrt{\eps} \bigg(\int_{\RR^3} |\bigg(\nabla, \frac{\partial}{\partial z} \bigg) U|^2 \bigg)^{1/2}.
\end{align}
By \eqref{stra4},  \eqref{stratil} and \eqref{stradep}, we conclude that the total stray field energy is bounded by:
$$\int_{\RR^3} |\bigg(\nabla, \frac{\partial}{\partial z} \bigg) U|^2 \leq  \frac{\pi+o(1)}{|\log
\eta|},$$
i.e., by \eqref{choiceUuu},
\be
\label{stratotal}
\frac 1 \eta \int_{\RR^2}
\left|\,|\nabla|^{-1/2}(\nabla\cdot m')\right|^2\, dx\leq \frac{2\pi+o(1)}{\eta |\log
\eta|}.
\ee
Finally, we estimate the last term of our energy given by the $m_3-$component.
For our configuration $m$, the only region in $\Omega$ where $m$ is not in-plane corresponds to the vortex core $B_\eps$. There we have
$$\frac{1}{\eps^2} \int_\Omega m_3^2 \, dx=\frac{1}{\eps^2} \int_{B_\eps} \cos^2 (\frac{\pi}{2\eps}r)\, dx=O(1).$$
Combining with \eqref{enextotal} and \eqref{stratotal}, the conclusion follows. Remark that the constructed configuration $m\in H^1(\Omega, S^2)$ is only continuous. By the density of $C^1(\Omega, S^2)$ vector fields satisfying \eqref{cond_tan} in the space of $H^1(\Omega, S^2)$ vector fields with \eqref{cond_tan} (for $C^{1,1}$ domains), one can smooth the configuration $m$ so that
the previous upper bound remains true.
\qed

\section{Appendix}

As mentioned in introduction, condition \eqref{cond_tan} is necessary for a configuration to have finite stray-field energy in our model. To simplify the notation, we prove the statement for the case where $\partial \Omega$ is a straight line:

\begin{pro}
\label{prop_app}
Let $\Omega=(-\infty, 0)\times \RR$ and $m'\in H^1(\Omega, \RR^2)$. With the convention $m':=m'{\bf 1}_\Omega$, then
$$\int_{\RR^2}
\left|\,|\nabla|^{-1/2}(\nabla\cdot m')\right|^2\, dx<\infty \quad \textrm{ implies that } \quad m_1(0, \cdot)=0 \textrm{ in } H^{1/2}(\RR).$$
\end{pro}

\proof{}
We will show that
\be
\label{claimi}
\int_{\RR} m_1(0, x_2) \f(x_2)\, dx_2=0 \textrm{ for every } \f\in C^\infty_c(\RR).\ee
(Here, $m_1(0, \cdot)$ represents the $H^{1/2}(\RR)-$trace on the vertical line $\{x_1=0\}$ of $m_1\in H^1(\Omega, \RR)$).
For a small $\eps>0$, let $\zeta_\eps$ mimic the normal component of a N\'eel wall transition on a scale $\eps$ with the size of the
core of order ${\eps^2}$:
$$\zeta_\eps(x_1)=\begin{cases} \frac{\log \frac{\eps^2}{(x_1^2+\eps^4)}}{\log \frac{1}{\eps^2}}&\quad \textrm{ if } |x_1|\leq \sqrt{\eps^2-\eps^4},\\
0 &\quad \textrm{ elsewhere}.
\end{cases}$$
We claim that \eqref{claimi} is equivalent to 
\be
\label{zer1}
\lim_{\eps\to 0}\int_{\RR^2} m_1(x_1, x_2) \f(x_2)\frac{d\ze}{dx_1}(x_1)\, dx_1 dx_2=0 \textrm{ for every } \f\in C^\infty_c(\RR).\ee
Indeed, we have:
\begin{align*}
&\left| \int_{\RR^2} m_1(x_1, x_2) \f(x_2)\frac{d\ze}{dx_1}(x_1)\, dx_1 dx_2-\int_{\RR} m_1(0, x_2) \f(x_2)\, dx_2\right|\\
\leq & \left| \int_{\RR^2} m_1(x_1, x_2) \f(x_2)\frac{d\ze}{dx_1}(x_1)\, dx_1 dx_2-\int_{-\infty}^0 \frac{d\ze}{dx_1}(x_1)\int_{\RR} m_1(0, x_2) \f(x_2)\, dx_2 dx_1\right|\\
= & \left| \int_{-\eps}^0  \frac{d\ze}{dx_1}(x_1)  \int_{\RR} \f(x_2) \bigg( \int_{x_1}^0 \frac{\partial m_1}{\partial x_1}(s, x_2) \, ds\bigg) \, dx_2 dx_1\right|\\
\leq & \int_{\RR} |\f(x_2)| \int_{-\eps}^0 \bigg|\frac{\partial m_1}{\partial x_1}(s, x_2)\bigg| \, ds dx_2\\
\leq & \sqrt{\eps} \int_{\RR} |\f(x_2)| \|\frac{\partial m_1}{\partial x_1}(\cdot, x_2)\|_{L^2(\RR_-)} \, dx_2\leq \sqrt{\eps} \|\f\|_{L^2(\RR)} \|\frac{\partial m_1}{\partial x_1}\|_{L^2(\Omega)}
\end{align*}
(where we used that $\ze$ is increasing on $\RR_-$ and $\int_{\RR_-} \frac{d\ze}{dx_1}(x_1)\, dx_1=1$).
In order to prove \eqref{zer1}, we set $\psi(x_1, x_2)=\ze(x_1) \f(x_2)$ and we write
$$\int_{\RR^2} m_1(x_1, x_2) \f(x_2)\frac{d\ze}{dx_1}(x_1)\, dx_1 dx_2=\int_{\RR^2} m'\cdot \nabla \psi \, dx_1 dx_2-\int_{\RR^2} m_2(x_1, x_2) \ze(x_1)\frac{d\f}{dx_2}(x_2)\, dx_1 dx_2.$$
Integrating by parts, we estimate the second term in the above RHS:
\begin{align*}
\left| \int_{\RR^2} m_2(x_1, x_2) \ze(x_1)\frac{d\f}{dx_2}(x_2)\, dx_1 dx_2\right|&=\left| \int_{-\eps}^0 \ze(x_1) \int_{\RR} \f(x_2)\frac{\partial m_2}{\partial x_2}(x_1, x_2)\, dx_2 dx_1\right|\\
&\leq  \|\f\|_{L^2(\RR)} \int_{-\eps}^0 \ze(x_1) \|\frac{\partial m_2}{\partial x_2}(x_1, \cdot)\|_{L^2(\RR)} \, dx_1\\
&\leq \sqrt{\eps}\|\f\|_{L^2(\RR)} \|\frac{\partial m_2}{\partial x_2}\|_{L^2(\Omega)}
\end{align*}
(since $\int_{-\eps}^0 \ze^2(x_1)\, dx_1\stackrel{x_1=\eps s}{\leq} \frac{C\eps}{|\log \eps|^2} \int_0^{\sqrt{1-\eps^2}} \log^2(s^2+\eps^2)\, ds=O(\frac{\eps}{|\log \eps|^2})\, $).
The first term in the above RHS is estimated by interpolation:
$$\left|\int_{\RR^2} m'\cdot \nabla \psi \right|=\left|\int_{\RR^2}\nabla\cdot m' \psi \right|\leq \int_{\RR^2}
\left|\,|\nabla|^{-1/2}(\nabla\cdot m')\right|^2 \int_{\RR^2}
\left|\,|\nabla|^{1/2}\psi\right|^2.$$
In order to conclude, we need to prove that $\|\psi\|_{\dot{H}^{1/2}(\RR^2)}\to 0$ as $\eps\to 0$.
For that, we use \eqref{esti_h12_ex} (valid in any dimension) for the following extension $V:\RR^3\to \RR$ of $\psi$ given by $V(x_1, x_2, z)=\psi(r, x_2)=\ze(r)\f(x_2)$ for every $(x_1, x_2, z)\in \RR^3$ and $r=\sqrt{x_1^2+z^2}$:
$$\big|\nabla V\big|^2+\big|\frac{\partial
V}{\partial z}\big|^2=\ze^2(r) \left| \frac{d\f}{dx_2}(x_2)\right|^2+\f^2(x_2)\left| \frac{d\ze}{dr}(r)\right|^2$$
and
\begin{align*}
\frac 1 \pi \int_{\RR^2}\left|\,|\nabla|^{1/2}\psi\right|^2\, dx&\leq \frac{1}{2\pi} \int_{\RR^3} \bigg(\big|\nabla V\big|^2+\big|\frac{\partial
V}{\partial z}\big|^2 \bigg)\, dxdz\\
&=  \|\frac{d\f}{d x_2}\|^2_{L^2(\RR)} \int_0^{\eps} r\ze^2(r)\, dr+\|\f\|^2_{L^2(\RR)} \int_0^{\eps} r\left| \frac{d\ze}{dr}(r)\right|^2\, dr\\
&\stackrel{\eqref{Straytran}}{\leq} C\bigg(\eps^2  \|\frac{d\f}{d x_2}\|^2_{L^2(\RR)}+\frac{1}{|\log \eps|} \|\f\|^2_{L^2(\RR)}\bigg)\to 0 \, \, \textrm{ as } \eps\to 0.
\end{align*}
\qed

\bigskip
\nd {\bf Acknowledgement}: 
R.I. \& F.O. acknowledge partial support of the German Science Foundation (DFG) through the Collaborative Research Center SFB 611. R.I. thanks hospitality to
Hausdorff Research Institute (Bonn) and Max Planck Institute (Leipzig).

\bibliographystyle{mystyle_fr}
\bibliography{references}

\end{document}